\documentclass{amsart}
\usepackage{graphicx}
\usepackage[all]{xy}
\usepackage{latexsym}
\usepackage{amssymb}
\usepackage{amsmath}
\usepackage{color}

\xyoption{arc}

 \newfam\cyrfam

  \font\tencyr=wncyr10

  \font\sevencyr=wncyr7

  \font\fivecyr=wncyr5

  \textfont\cyrfam=\tencyr \scriptfont\cyrfam=\sevencyr

    \scriptscriptfont\cyrfam=\fivecyr

  \newfam\cyifam

  \font\tencyi=wncyi10

  \font\sevencyi=wncyi7

  \font\fivecyi=wncyi5

  \textfont\cyifam=\tencyi \scriptfont\cyifam=\sevencyi

    \scriptscriptfont\cyifam=\fivecyi

\def\id{{\mbox{1 \hskip -7pt 1}}}
\newcommand{\sgn}{{\mathit s  \mathit g\mathit  n}}
 \newcommand{\lon}{\longrightarrow}
 \newcommand{\bu}{\bullet}
 
 \newcommand{\rar}{\rightarrow}
 \newcommand{\hook}{\hookrightarrow}

\newcommand{\p}{{\partial}}

\newcommand{\Der}{\mathrm{Der}}

\newcommand{\tw}{\mathsf{tw}}
\newcommand{\RGra}{\mathcal{R} \mathcal{G} ra}

 \newcommand{\Z}{{\mathbb Z}}
 \newcommand{\bS}{{\mathbb S}}

 \newcommand{\K}{{\mathbb K}}

 \newcommand{\ot}{\otimes}

\newcommand{\Def}{\mathsf{Def}}

\newcommand{\GC}{\mathsf{GC}}
\newcommand{\OGC}{\mathsf{OGC}}

\newcommand{\grt}{\fg\fr\ft}
 \newcommand{\Beq}{\begin{equation}}
 \newcommand{\Eeq}{\end{equation}}
 \newcommand{\Beqr}{\begin{eqnarray}}
 \newcommand{\Eeqr}{\end{eqnarray}}
 \newcommand{\Beqrn}{\begin{eqnarray*}}
 \newcommand{\Eeqrn}{\end{eqnarray*}}
 \newcommand{\Ba}{\begin{array}}
 \newcommand{\Ea}{\end{array}}
 \newcommand{\Bi}{\begin{itemize}}
 \newcommand{\Ei}{\end{itemize}}
 \newcommand{\Bc}{\begin{center}}
 \newcommand{\Ec}{\end{center}}

 \newcommand{\fg}{{\mathfrak g}}

\newcommand{\fr}{{\mathfrak r}}

\newcommand{\ft}{{\mathfrak t}}

 \newcommand{\cA}{{\mathcal A}}

 \newcommand{\cE}{{\mathcal E}}
 \newcommand{\cF}{{\mathcal F}}
 \newcommand{\cG}{{\mathcal G}}

 \newcommand{\cM}{{\mathcal M}}

 \newcommand{\cR}{{\mathcal R}}

 \newcommand{\ga}{\gamma}
 
 \newcommand{\Ga}{\Gamma}

 \newcommand{\la}{\lambda}

 \newcommand{\Ker}{{\mathsf K \mathsf e \mathsf r}\, }
 
 \newcommand{\Hom}{{\mathrm H\mathrm o\mathrm m}}
 
 \newcommand{\sip}{\smallskip}
 \newcommand{\bip}{\bigskip}
 \newcommand{\mip}{\vspace{2.5mm}}

\newcommand{\LB}{\mathcal{L}\mathit{ieb}}

\newcommand{\HoLB}{\mathcal{H}\mathit{olieb}}

\newcommand{\RGC}{\mathsf{RGC}}

\newcommand{\Holie}{\mathcal{H} \mathit{olie}}
\newcommand{\Lie}{\mathcal{L} \mathit{ie}}

\newcommand{\Gra}{\mathcal{G}ra}
\newcommand{\DRGra}{\mathcal{RG}ra^\uparrow}
\newcommand{\RGraphs}{\mathcal{RG}raphs}

\newcommand{\ORGra}{\mathcal{RG}ra^{or}}
\newcommand{\ORGraphs}{\mathcal{RG}raphs^{or}}
\newcommand{\OORGraphs}{O\mathcal{RG}raphs}

 \newcommand{\bbu}{\mbox{\resizebox{2.8mm}{!}{$\bullet$}}}

\newcommand{\ORGC}{\mathsf{ORGC}}
\newcommand{\PCY}{\mathcal{P}re\mathcal{CY}^3}
\newcommand{\PCYf}{\mathcal{P}re\mathcal{CY}}

\newcommand{\rgc}{\mathfrak{rgc}}
\newcommand{\orgc}{\mathfrak{orgc}}

\theoremstyle{plain}
\swapnumbers

\newtheorem{prop-def}[theorem]{Proposition-definition}

\newtheorem{f-theorem}{Formality Theorem}[section]
\newtheorem{main-theorem}{Main~Theorem}[section]
\newtheorem{section-theorem}{Theorem}[section]

\theoremstyle{definition}

\usepackage{tikz}
\usetikzlibrary{matrix}
\usetikzlibrary{shapes}
\usetikzlibrary{arrows}
\usetikzlibrary{calc,3d}
\usetikzlibrary{decorations,decorations.pathmorphing,decorations.pathreplacing,decorations.markings}
\usetikzlibrary{through}

\tikzset{ext/.style={circle, draw,inner sep=1pt},int/.style={circle,draw,fill,inner sep=1.4pt},nil/.style={inner sep=1pt}}
\tikzset{cy/.style={circle,draw,fill,inner sep=2pt},scy/.style={circle,draw,inner sep=2pt},scyx/.style={draw,cross out,inner sep=2pt},scyt/.style={draw,regular polygon,regular polygon sides=3,inner sep=0.95pt}}
\tikzset{exte/.style={circle, draw,inner sep=3pt},inte/.style={circle,draw,fill,inner sep=3pt}}
\tikzset{diagram/.style={matrix of math nodes, row sep=3em, column sep=2.5em, text height=1.5ex, text depth=0.25ex}}
\tikzset{diagram2/.style={matrix of math nodes, row sep=0.5em, column sep=0.5em, text height=1.5ex, text depth=0.25ex}}

 \tikzset{
  rightblack/.style={
    decoration={markings,mark=at position .8 with {\arrow[scale=1.2,black]{latex}}},
    postaction={decorate},
    shorten >=0.4pt}}
\tikzset{
  leftblack/.style={
    decoration={markings,mark=at position .55 with {\arrowreversed[scale=1.2,black]{latex}}},
    postaction={decorate},
    shorten >=0.4pt}}

\begin{document}

 \sloppy

 \newenvironment{proo}{\begin{trivlist} \item{\sc {Proof.}}}
  {\hfill $\square$ \end{trivlist}}

\long\def\symbolfootnote[#1]#2{\begingroup%
\def\thefootnote{\fnsymbol{footnote}}\footnote[#1]{#2}\endgroup}

\title{A complex of ribbon quivers
 and $\cM_{g,m}$}

\author{Sergei A.\ Merkulov}
\email{smerkulov25@gmail.com}

\begin{abstract}  For any integer $d\in \Z$   
we introduce a  complex $\ORGC_{d}^{(g,m)}$ spanned 
by genus $g$  ribbon quivers with $m$ marked boundaries and prove that its cohomology 
computes (up to a degree shift)  the  compactly supported cohomology of the moduli 
space $\cM_{g,m}$ of genus $g$ algebraic curves  with $m$  marked points.

\sip

We show that the totality of complexes 
$$
\orgc_{d}= \prod_{g\geq 1} \ORGC_{d}^{(g,1)} \stackrel{qis}{\simeq}  
\prod_{g\geq 1} H_c^{\bu-1+2g(d-1)}(\cM_{g,1})
$$ 
has a natural dg Lie algebra structure
 which controls the deformation theory of the dg properad $\PCY_d$ governing a certain class of
 (possibly, infinite-dimensional) degree $d$ pre-Calabi-Yau algebras. 
 This result implies
 in particular that for $d\leq 2$ the zero-th cohomology group of the  
 complex $\Der(\PCY_d)$ of the preserving boundaries derivations is one-dimensional, while for $d=2$ the cohomology group $H^1(\Der(\PCY_2))$ 
  contains a subspace isomorphic to the Grothendieck-Teichm\"uller Lie algebra.


\end{abstract}

 \maketitle
\markboth{}{}

{\Large
\section{\bf Introduction}
}

\
\subsection{A new combinatorial model for  the cohomology of $\cM_{g,m}$}
For any integer $d$ there is a cochain complex $\RGC_d$  generated by ribbon graphs with no
fixed directions on their edges\footnote{More precisely, for $d$ even each edge is  undirected,
while for $d$ odd each edge is equipped with an ordering of its two half-edges up to a flip and
multiplication by $-1$. In both cases edges are assigned the cohomological degree $1-d$,
see \S 2 for full details.} which computes the totality
of compactly supported cohomology groups \cite{Pe} 
$$
H^\bu(\RGC_d)=\prod_{g\geq 0, m\geq 1\atop 2g+m\geq 3} H^{\bu +(d-1)m+d(2g-1)}_c(\cM_{g,m})
$$
of the moduli spaces $\cM_{g,m}$ of genus $g$ algebraic curves with $m$ marked points.

\sip

For any $d\in \Z$ we introduce in this paper a new cochain complex $\ORGC_d$ generated by
ribbon graphs $\Ga$  whose edges are assigned the same cohomological degree $1-d$ as in the
case of $\RGC_d$ but, in contrast to the latter, the directions of edges of $\Ga$  are
strictly fixed and, most importantly,   {\it no}\, oriented closed paths of directed edges are allowed in $\Ga$; such ribbon graphs are called {\it ribbon quivers}. By analogy to
the ``commutative" (i.e.\ non-ribbon) case studied in \cite{Wi2} one may also call such ribbon graphs
{\it oriented}: they have a well-defined directed flow on all edges. We prove in this
paper that this complex of ribbon quivers is quasi-isomorphic
to R.\ Penner's ribbon graph complex.

 \subsubsection{\bf Theorem}\label{1: Main theorem 1} {\it For any $d\in \Z$ one has}
\Beq\label{1: H(RGC)=H(ORGC)}
 H^{\bu}(\RGC_d)=H^\bu(\ORGC_{d+1}).
\Eeq

\sip

Notice the degree shift $d+1$ in the r.h.s.\ of this equation. Thus the complex of ribbon quivers
gives us a new combinatorial model for the totality of compactly supported cohomology groups of the moduli spaces
$\cM_{g,m}$. This result is a ribbon analogue of the similar isomorphism \cite{Wi2,Z,Me3, MWW}
$$
H^\bu(\GC_d)=H^\bu(\OGC_{d+1})
$$
which holds true for the famous M.\ Kontsevich's  \cite{Ko2} ``commutative" graph complex 
$\GC_d$ and its oriented version $\OGC_{d}$. The graph complex $\GC_d$  found 
many important applications 
in various areas of mathematics and mathematical physics (see e.g.\ the expositions \cite{Ko3,Wi3} 
and references cited there). The complex $\OGC_{d+1}$ --- the ``commutative" analogue of 
$\ORGC_{d+1}$ --- was also studied earlier \cite{Wi2,Z,Me2, MWW} as it admits 
applications in the homotopy theory of (involutive) Lie bialgebras and their 
quantizations, in string topology and some other areas of modern research. 

\subsection{Lie algebra structures on ribbon graph complexes}
  Both complexes $\RGC_d$ and $\ORGC_{d+1}$ split into direct products of subcomplexes,
$$
\RGC_d=\prod_{g\geq 0, m\geq 1\atop 2g+m\geq 3} \RGC_d^{(g,m)}, \ \ \ \
\ORGC_{d+1}=\prod_{g\geq 0, m\geq 1\atop 2g+m\geq 3} \ORGC_{d+1}^{(g,m)},
$$
generated by  genus $g$ ribbon graphs with precisely $m$ marked boundaries, and
the above isomorphism (\ref{1: H(RGC)=H(ORGC)}) respects these decompositions,
$$
 H^\bu(\RGC_d^{(g,m)})=H^\bu(\ORGC_{d+1}^{(g,m)})\simeq  H^{\bu +(d-1)m+d(2g-1)}_c(\cM_{g,m}).
$$
It was shown in \cite{MW1} that the totality of (degree shifted) subcomplexes of $\RGC_d$ generated by 
ribbon graphs with precisely one
marked point,
$$
\rgc_d:= \prod_{g\geq 1}\RGC_d^{(g,1)} \stackrel{qis}{\simeq}\prod_{g\geq 1} H^{\bu -1+ 2dg}_c(\cM_{g,1})
$$
 has a natural dg Lie algebra structure given by substitution of the unique boundary of one ribbon graph into vertices of another ribbon graph (see \S 2.2 below). Moreover,  Penner's complex $\RGC_d$ is 
 a $\rgc_d$-module. The same is true for the oriented ribbon graph complex: 
 the totality
  $$
\orgc_d:= \prod_{g\geq 1}\ORGC_d^{(g,1)}.
$$
is a dg Lie algebra, and $\ORGC_d$ is an $\orgc_d$-module. An important for applications point
is that the isomorphism (\ref{1: H(RGC)=H(ORGC)}) respects these extra structures. 

 \subsubsection{\bf Theorem}\label{1: Main theorem 2}
  {\it  The dg Lie algebras $\rgc_d$ and $\orgc_{d+1}$ are $\Lie_\infty$ quasi-isomorphic for any $d\in \Z$.}

 \sip
  
  This result computes the cohomology of $\orgc_d$ in terms of the compactly supported cohomology groups of the moduli spaces $\cM_{g,1}$,
  $$
H^\bu(\orgc_{d})=  
\prod_{g\geq 1} H_c^{\bu-1+2g(d-1)}(\cM_{g,1}).
$$

  \subsection{An application to the theory of pre-Calabi-Yau structures}
Pre-Calabi-Yau algebra structures in a dg vector space $A$ are defined as Maurer-Cartan
 elements $\pi$ of the so called higher Hochschield graded Lie algebra \cite{IK,IKV,KTV}
$$
{C}_{[d]}(A) := \prod_{k\geq 1}
\left( \bigoplus_{n_1,\ldots, n_k\geq 0}
\Hom\left(\bigotimes_{i=1}^{k}
(A[1])^{\ot n_i}, (A[2-d])^{\ot k}\right)_{\Z_k}\right).
$$
 There is a dg free properad  $\PCYf_d$ which is uniquely characterized  the 
 following property: one has a one-to-one correspondence between representations of 
 $\PCYf_d$ in a dg vector space $A$ and pre-Calabi-Yau algebra  structures $\pi$ in $A$; it has been explicitly described in \cite{LV}. 
 We study in this paper the universal deformation theory 
 of a class of pre-Calabi-Yau algebras which  are defined as  MC elements of the 
 following Lie subalgebra  of ${C}_{[d]}(A)$,
  $$
{C}_{[d]}^3(A) := \prod_{k\geq 1}
\left( \bigoplus_{n_1,\ldots, n_k\geq 0\atop {n_1+...+n_k\geq 1\atop n_1+..+n_k+k\geq 3}}
\Hom\left(\bigotimes_{i=1}^{k}
(A[1])^{\ot n_i}, (A[2-d])^{\ot k}\right)_{\Z_k}\right).
$$
  The superscript $3$  emphasizes the imposed ``at least trivalency" condition  
  $ n_1+..+n_k+k\geq 3$ which guarantees, in particular, that the $k=1$ part of $\pi$ is a genuine 
  $A_\infty$ algebra structure on $A$. These structures are controlled by a dg free
   properad $\PCY_d$ which has been first explicitly described in \cite{Q} as the cobar construction of the cooperad Koszul dual to the properad $\mathsf{BIB}_d$ of so called balanced infinitesimal bialgebras; it was also shown in \cite{Q} that $\PCY_2$ admits an epimorphism,
  $$
  \PCY_2 \lon \mathsf{DPoiss}_\infty,
  $$
 into the dg properad  $\mathsf{DPoiss}_\infty$  governing strongly homotopy double Poisson algebras \cite{LV}. 
 The properad $\PCY_d$ is a kind of non-commutative version of the dg properad $\HoLB_d$  which  controls dergee $d$ shifted  strongly homotopy Lie bialgebras  (whose deformation theory was studied in \cite{MW2}).

 \sip

 In our paper  the properad $\PCYf_d$ and its quotient properad $\PCY_d$ emerge as dg properads of ribbon quivers 
 {\it with hairs}\, (see \S 5 below) so that it makes sense talking about boundaries of its elements.  We study in this 
 paper the genus completed dg Lie algebra $\Der(\PCY_d)$ of preserving the boundaries 
 derivations (see \S {\ref{5: remaik on full hairs}} for the precise definition) and compute its cohomology in terms of $H^\bu_c(\cM_{g,1})$. More precisely
 we prove the following statement.

  \subsubsection{\bf Theorem}\label{1: Theorem on H(Der)}  {\it There is an explicit
 morphism of dg Lie algebras
  $$
  \orgc_d \lon \Der(\PCY_d)
  $$ 
  which is a quasi-isomorphism up to one rescaling cohomology class. Hence
  $$
  H^\bu(\Der(\PCY_d))=\prod_{g\geq 1} H_c^{\bu-1+2g(d-1)}(\cM_{g,1}) \oplus \K[0]
  $$
   }
  
  \sip
  
This result implies, in particular, that $H^0(\Der(\PCY_d))=\K$ for any $d\leq 2$,
i.e.\ that the dg properad $\PCY_{d}$ has  no homotopy non-trivial automorphisms (except rescalings) for $d\leq 2$
which preserve the number of boundaries (but not the genus). On the other hand
for $d=2$ one concludes that $H^{1}(\Der(\PCY_2))$ contains a subspace isomorphic 
to the Grothendieck-Teichm\"uller Lie algebra $\grt$ (see \S 5.6 below).

\sip

As another  application of the above results we show in \S 5.8 below that the above mentioned properad $\mathsf{BIB}_d$ of degree $d$ balanced infinitesimal bialgebras which was introduced and studied  by Alexandre Quesney  
in \cite{Q}, is not Koszul.   

  \subsection{Some notation} We work over a field $\K$ of characteristic zero.
 The set $\{1,2, \ldots, n\}$ is abbreviated to $[n]$;  its group of automorphisms is
denoted by $\bS_n$; the trivial (resp., the sign) one-dimensional representation of
 $\bS_n$ is denoted by $\id_n$ (resp.,  $\sgn_n$). The cardinality of a finite set $S$ is 
 denoted by $\# S$ while its linear span over a
field $\K$ by $\K\left\langle S\right\rangle$.
If $V=\oplus_{i\in \Z} V^i$ is a graded vector space, then
$V[k]$ stands for the graded vector space with $V[k]^i:=V^{i+k}$. For 
$v\in V^i$ we set $|v|:=i$.

\mip

{\bf Acknowledgement}. It is a great pleasure to thank
Alexandre Quesney for very usefull comments and discussions.

\bip

{\Large
\section{\bf Props and dg Lie algebras of ribbon graphs}
}

\mip

\subsection{Reminder on the operad of (curved) $\Lie_d$-algebras} 
Recall that the operad of degree shifted Lie algebras is defined, for any integer  $d\in \Z$, 
as the quotient,
$$
\Lie_{d}:=\cF ree\langle E\rangle/\langle\cR\rangle,
$$
of the free operad generated by an  $\bS$-module $E=\{E(n)\}_{n\geq 2}$ with
 all $E(n)=0$ except
$$
E(2):= \sgn_2^{|d|}[d-1]=\mbox{span}\left\langle
\Ba{c}\begin{xy}
 <0mm,0.66mm>*{};<0mm,3mm>*{}**@{-},
 <0.39mm,-0.39mm>*{};<2.2mm,-2.2mm>*{}**@{-},
 <-0.35mm,-0.35mm>*{};<-2.2mm,-2.2mm>*{}**@{-},
 <0mm,0mm>*{\bu};<0mm,0mm>*{}**@{},
   <0.39mm,-0.39mm>*{};<2.9mm,-4mm>*{^{_2}}**@{},
   <-0.35mm,-0.35mm>*{};<-2.8mm,-4mm>*{^{_1}}**@{},
\end{xy}\Ea
=(-1)^{d}
\Ba{c}\begin{xy}
 <0mm,0.66mm>*{};<0mm,3mm>*{}**@{-},
 <0.39mm,-0.39mm>*{};<2.2mm,-2.2mm>*{}**@{-},
 <-0.35mm,-0.35mm>*{};<-2.2mm,-2.2mm>*{}**@{-},
 <0mm,0mm>*{\bu};<0mm,0mm>*{}**@{},
   <0.39mm,-0.39mm>*{};<2.9mm,-4mm>*{^{_1}}**@{},
   <-0.35mm,-0.35mm>*{};<-2.8mm,-4mm>*{^{_2}}**@{},
\end{xy}\Ea
\right\rangle,
$$
modulo the ideal generated by the following relation
\Beq\label{2: Lie operad Jacobi relation}
 \Ba{c}\resizebox{9mm}{!}{ \begin{xy}
 <0mm,0mm>*{\bu};<0mm,0mm>*{}**@{},
 <0mm,0.69mm>*{};<0mm,3.0mm>*{}**@{-},
 <0.39mm,-0.39mm>*{};<2.4mm,-2.4mm>*{}**@{-},
 <-0.35mm,-0.35mm>*{};<-1.9mm,-1.9mm>*{}**@{-},
 <-2.4mm,-2.4mm>*{\bu};<-2.4mm,-2.4mm>*{}**@{},
 <-2.0mm,-2.8mm>*{};<0mm,-4.9mm>*{}**@{-},
 <-2.8mm,-2.9mm>*{};<-4.7mm,-4.9mm>*{}**@{-},
    <0.39mm,-0.39mm>*{};<3.3mm,-4.0mm>*{^3}**@{},
    <-2.0mm,-2.8mm>*{};<0.5mm,-6.7mm>*{^2}**@{},
    <-2.8mm,-2.9mm>*{};<-5.2mm,-6.7mm>*{^1}**@{},
 \end{xy}}\Ea
 +
\Ba{c}\resizebox{9mm}{!}{ \begin{xy}
 <0mm,0mm>*{\bu};<0mm,0mm>*{}**@{},
 <0mm,0.69mm>*{};<0mm,3.0mm>*{}**@{-},
 <0.39mm,-0.39mm>*{};<2.4mm,-2.4mm>*{}**@{-},
 <-0.35mm,-0.35mm>*{};<-1.9mm,-1.9mm>*{}**@{-},
 <-2.4mm,-2.4mm>*{\bu};<-2.4mm,-2.4mm>*{}**@{},
 <-2.0mm,-2.8mm>*{};<0mm,-4.9mm>*{}**@{-},
 <-2.8mm,-2.9mm>*{};<-4.7mm,-4.9mm>*{}**@{-},
    <0.39mm,-0.39mm>*{};<3.3mm,-4.0mm>*{^2}**@{},
    <-2.0mm,-2.8mm>*{};<0.5mm,-6.7mm>*{^1}**@{},
    <-2.8mm,-2.9mm>*{};<-5.2mm,-6.7mm>*{^3}**@{},
 \end{xy}}\Ea
 +
\Ba{c}\resizebox{9mm}{!}{ \begin{xy}
 <0mm,0mm>*{\bu};<0mm,0mm>*{}**@{},
 <0mm,0.69mm>*{};<0mm,3.0mm>*{}**@{-},
 <0.39mm,-0.39mm>*{};<2.4mm,-2.4mm>*{}**@{-},
 <-0.35mm,-0.35mm>*{};<-1.9mm,-1.9mm>*{}**@{-},
 <-2.4mm,-2.4mm>*{\bu};<-2.4mm,-2.4mm>*{}**@{},
 <-2.0mm,-2.8mm>*{};<0mm,-4.9mm>*{}**@{-},
 <-2.8mm,-2.9mm>*{};<-4.7mm,-4.9mm>*{}**@{-},
    <0.39mm,-0.39mm>*{};<3.3mm,-4.0mm>*{^1}**@{},
    <-2.0mm,-2.8mm>*{};<0.5mm,-6.7mm>*{^3}**@{},
    <-2.8mm,-2.9mm>*{};<-5.2mm,-6.7mm>*{^2}**@{},
 \end{xy}}\Ea
 =0.
\Eeq
Its minimal resolution $\Holie_d$ is a dg free operad whose (skew)symmetric generators,
\Beq\label{2: Lie_inf corolla}
\Ba{c}\resizebox{22mm}{!}{ \xy
(1,-5)*{\ldots},
(-13,-7)*{_1},
(-8,-7)*{_2},
(-3,-7)*{_3},
(7,-7)*{_{n-1}},
(13,-7)*{_n},
 (0,0)*{\bu}="a",
(0,5)*{}="0",
(-12,-5)*{}="b_1",
(-8,-5)*{}="b_2",
(-3,-5)*{}="b_3",
(8,-5)*{}="b_4",
(12,-5)*{}="b_5",
\ar @{-} "a";"0" <0pt>
\ar @{-} "a";"b_2" <0pt>
\ar @{-} "a";"b_3" <0pt>
\ar @{-} "a";"b_1" <0pt>
\ar @{-} "a";"b_4" <0pt>
\ar @{-} "a";"b_5" <0pt>
\endxy}\Ea
=(-1)^d
\Ba{c}\resizebox{23mm}{!}{\xy
(1,-6)*{\ldots},
(-13,-7)*{_{\sigma(1)}},
(-6.7,-7)*{_{\sigma(2)}},
(13,-7)*{_{\sigma(n)}},
 (0,0)*{\bu}="a",
(0,5)*{}="0",
(-12,-5)*{}="b_1",
(-8,-5)*{}="b_2",
(-3,-5)*{}="b_3",
(8,-5)*{}="b_4",
(12,-5)*{}="b_5",
\ar @{-} "a";"0" <0pt>
\ar @{-} "a";"b_2" <0pt>
\ar @{-} "a";"b_3" <0pt>
\ar @{-} "a";"b_1" <0pt>
\ar @{-} "a";"b_4" <0pt>
\ar @{-} "a";"b_5" <0pt>
\endxy}\Ea,
\ \ \ \forall \sigma\in \bS_n,\ n\geq2,
\Eeq
have degrees $1+d-nd$.
 The differential in $\Holie_d$ is given by
\Beq\label{2: d in Lie_infty}
\delta\hspace{-3mm}
\Ba{c}\resizebox{21mm}{!}{\xy
(1,-5)*{\ldots},
(-13,-7)*{_1},
(-8,-7)*{_2},
(-3,-7)*{_3},
(7,-7)*{_{n-1}},
(13,-7)*{_n},
 (0,0)*{\bu}="a",
(0,5)*{}="0",
(-12,-5)*{}="b_1",
(-8,-5)*{}="b_2",
(-3,-5)*{}="b_3",
(8,-5)*{}="b_4",
(12,-5)*{}="b_5",
\ar @{-} "a";"0" <0pt>
\ar @{-} "a";"b_2" <0pt>
\ar @{-} "a";"b_3" <0pt>
\ar @{-} "a";"b_1" <0pt>
\ar @{-} "a";"b_4" <0pt>
\ar @{-} "a";"b_5" <0pt>
\endxy}\Ea
=
\sum_{A\varsubsetneq [n]\atop
\# A\geq 2}\pm
%
%
\Ba{c}\resizebox{19mm}{!}{\begin{xy}
<10mm,0mm>*{\bu},
<10mm,0.8mm>*{};<10mm,5mm>*{}**@{-},
<0mm,-10mm>*{...},
<14mm,-5mm>*{\ldots},
<13mm,-7mm>*{\underbrace{\ \ \ \ \ \ \ \ \ \ \ \ \  }},
<14mm,-10mm>*{_{[n]\setminus A}};
<10.3mm,0.1mm>*{};<20mm,-5mm>*{}**@{-},
<9.7mm,-0.5mm>*{};<6mm,-5mm>*{}**@{-},
<9.9mm,-0.5mm>*{};<10mm,-5mm>*{}**@{-},
<9.6mm,0.1mm>*{};<0mm,-4.4mm>*{}**@{-},
<0mm,-5mm>*{\bu};
<-5mm,-10mm>*{}**@{-},
<-2.7mm,-10mm>*{}**@{-},
<2.7mm,-10mm>*{}**@{-},
<5mm,-10mm>*{}**@{-},
<0mm,-12mm>*{\underbrace{\ \ \ \ \ \ \ \ \ \ }},
<0mm,-15mm>*{_{A}}.
\end{xy}}
\Ea
\Eeq
If $d$ is even, all the signs above are equal to $-1$.
The case $d=1$ corresponds to the usual operad of strongly homotopy Lie algebras 
which is often denoted by $\Lie_\infty$.

\sip

The dg operad of {\em curved}\, $c\Holie_d$-algebras is, by definition, a dg free operad 
generated by (skew)symmetric $n$-corollas  (\ref{2: Lie_inf corolla}) of degrees $1+d-nd$ for any 
$n\geq 0$ (rather than for any  $n\geq 2$ as in the case of $\Holie_d$). The differential in 
$c\Holie_d$ is given formally by the same splitting formula (\ref{2: d in Lie_infty}) with 
the condition $\#A\geq 2$ replaced by $\#A \geq 0$.


\subsection{Ribbon graphs}
A {\em ribbon graph}\, $\Ga$ is, by definition, a triple
$(H(\Ga), \sigma_1, \sigma_0)$ where
\Bi
\item[(i)]  $H(\Ga)$ is a finite set  called the
set of {\em half edges};
\item[(ii)] $\sigma_1$ a fixed point free involution $\sigma_1: H(\Ga)\rar H(\Ga)$ whose
set of orbits, $E(\Ga):=H(\Ga)/\sigma_1$, is called the set of  {\em edges};
\item[(iii)] $\sigma_0$ is an arbitrary  permutation $\sigma_0: H(\Ga)\rar H(\Ga)$. The set of orbits, $V(\Ga):=H(\Ga)/\sigma_0$,
is called the set of {\em vertices}\, of the ribbon graph. If we represent $\sigma_0$ as a product of cycles, then each cycle of $\sigma_0$ corresponds to a vertex $v\in V(\Ga)$ and the
preimage $H(v):=p^{-1}(v)\subset H(\Ga)$ under the projection
$$
p: H(\Ga) \rar V(\Ga)
$$
is called the set of half edges attached to $v$; this set comes equipped with a natural cyclic ordering.
\Ei
The orbits of the permutation $\sigma_\infty:= \sigma_0^{-1}\circ \sigma_1$
are called {\em boundaries}\, of the ribbon graph $\Ga$. The set of boundaries of $\Ga$ is
 denoted by $B(\Ga)$. The number 
$$
 g= 1+\frac{1}{2}\left(\# E(\Ga) - \# V(\Ga)- \# B(\Ga)\right)
 $$
is called the {\em genus}\, of a ribbon graph $\Ga$. 

\sip

There is an obvious geometric representation of a genus $g$ ribbon graph $\Ga$ as   a 2-dimensional 
genus $g$ topological surface with boundaries obtained by representing each vertex $v$ as a disk\footnote{In most of our pictures below we represent that
disk in a much reduced form as $\resizebox{2mm}{!}{\xy
 (10,1)*+{_i}*\frm{o}="B";
\endxy}$    or $\circ$ or $\bu$.}, 
\Beq\label{2: vertex v as a disk picture}
v= \Ba{c}\resizebox{20mm}{!}{\xy
 (0,0)*{
\xycircle(6,6){.}};
(-6,0)*{}="1";
(-16,0)*{}="1'";
(4.2,4.2)*{}="2";
(11,11)*{}="2'";
(4.2,-4.2)*{}="3";
(11,-11)*{}="3'";
(6,0)*{}="4";
(16,0)*{}="4'";
(-4.2,4.2)*{}="5";
(-11,11)*{}="5'";
(-4.2,-4.2)*{}="6";
(-12,-12)*{}="6'";
(-0.3,6)*{}="a";
(-0.4,6)*{}="b";
\ar @{-} "1";"1'" <0pt>
\ar @{-} "2";"2'" <0pt>
\ar @{-} "3";"3'" <0pt>
\ar @{-} "4";"4'" <0pt>
\ar @{-} "5";"5'" <0pt>
\ar @{->} "a";"b" <0pt>
\ar @{-} "6";"6'" <0pt>
\endxy}\Ea
\Eeq
and then thickening  each half-edge into a strip. Note that half-edges attached to a vertex  split the corresponding dashed circle
into a union of open dashed intervals,
$$
C(v)=\coprod c_i, \ \ \ 
v= \Ba{c}\resizebox{21mm}{!}{\xy
(8,2.5)*{{c_1}};
(0,8.8)*{{c_2}};
(-8,-2.5)*{{c_i}};
 (0,0)*{\xycircle(6,6){.}};
(-6,0)*{}="1";
(-16,0)*{}="1'";
(4.2,4.2)*{}="2";
(11,11)*{}="2'";
(4.2,-4.2)*{}="3";
(11,-11)*{}="3'";
(6,0)*{}="4";
(16,0)*{}="4'";
(-4.2,4.2)*{}="5";
(-11,11)*{}="5'";
(-4.2,-4.2)*{}="6";
(-11,-11)*{}="6'";
(-0.3,6)*{}="a";
(-0.4,6)*{}="b";
\ar @{-} "1";"1'" <0pt>
\ar @{-} "2";"2'" <0pt>
\ar @{-} "3";"3'" <0pt>
\ar @{-} "4";"4'" <0pt>
\ar @{-} "5";"5'" <0pt>
\ar @{->} "a";"b" <0pt>
\ar @{-} "6";"6'" <0pt>
\endxy}\Ea
\ \ \simeq \ \ 
\Ba{c}\resizebox{21mm}{!}{\xy
(6,2.1)*{{c_1}};
(0,3.8)*{{c_2}};
(-6,-2.5)*{{c_i}};
 (0,0)*{\bbu}="0";
(-6,0)*{}="1";
(-16,0)*{}="1'";
(4.2,4.2)*{}="2";
(11,11)*{}="2'";
(4.2,-4.2)*{}="3";
(11,-11)*{}="3'";
(6,0)*{}="4";
(16,0)*{}="4'";
(-4.2,4.2)*{}="5";
(-11,11)*{}="5'";
(-4.2,-4.2)*{}="6";
(-11,-11)*{}="6'";
\ar @{-} "0";"1'" <0pt>
\ar @{-} "0";"2'" <0pt>
\ar @{-} "0";"3'" <0pt>
\ar @{-} "0";"4'" <0pt>
\ar @{-} "0";"5'" <0pt>
\ar @{-} "0";"6'" <0pt>
\endxy}\Ea,
$$
 and we call each such interval a {\em corner}\, of the
ribbon graph $\Ga$ under consideration. Let us denote the set of all corners of $\Ga$ by
 $C(\Ga)$. There are two natural partitions of the set $C(\Ga)$,
 $$
 C(\Ga)=\coprod_{v\in V(\Ga)} C(v),\ \ \ \ \ \
 C(\Ga)=\coprod_{b\in B(\Ga)} C(b),\ \ \
 $$
 where each subset $C(v)$ or $C(b)$ comes equipped with a canonically induced cyclic ordering.
 For each boundary $b\in B(\Ga)$, the associated  cyclically ordered  set  $C(b)=\coprod_{k}c_k$ can be represented as a anti-clockwise oriented topological circle of the form
\Beq\label{2: boundary b picture}
 b=\Ba{c}\resizebox{18mm}{!}{\xy
 (14,0.5)*{{c_k}};
(-3,0.5)*{{c_j}};
(6,-8)*{\circlearrowleft},
 (0,0)*{}="a1",
(10,0)*{}="a2",
(13,-3)*{}="a3",
(13,-13)*{}="a4",
(10,-16)*{}="a5",
(0,-16)*{}="a6",
(-3,-13)*{}="a7",
(-3,-3)*{}="a8",
\ar @{-} "a1";"a2" <0pt>
\ar @{.} "a2";"a3" <0pt>
\ar @{-} "a3";"a4" <0pt>
\ar @{.} "a4";"a5" <0pt>
\ar @{-} "a5";"a6" <0pt>
\ar @{.} "a6";"a7" <0pt>
\ar @{-} "a7";"a8" <0pt>
\ar @{.} "a8";"a1" <0pt>
\endxy}\vspace{2mm}\Ea
\Eeq
where dashed intervals represent (always different) corners of $\Ga$ while
solid intervals
represent (not necessarily different) edges of $\Ga$ which belong to $b$ (cf.\ \cite{MW1}).

\subsection{A properad of ribbon graphs $\RGra_d$ and its twisted version}
Let $\RGra_d=\{\RGra_d(m,n)\}_{m,n\geq 1}$ be a properad of {\em connected}\, ribbon graphs introduced\footnote{More precisely, the symbol $\RGra_d$ stands in 
\cite{MW1} for the {\em prop}\, generated by not necessarily connected ribbon graphs; 
in this paper we work solely with connected
graphs and hence use the symbol $\RGra_d$ for the sub-properad of the latter which is generated by connected ribbon graphs.} in \S 4 of \cite{MW1}. The $\bS_m^{op}\times\bS_n$-module $\RGra_d(m,n)$
is generated by ribbon graphs $\Ga$ with $n$ labelled vertices and $m$ labelled
(by integers $\bar{1}, \bar{2},\ldots$) boundaries, e.g.
$$
\Ba{c} \resizebox{11mm}{!}{\xy
 (5,3)*{^{\bar{1}}};
 (10,1)*+{_1}*\frm{o}="B";
 (0,1)*+{_2}*\frm{o}="A";
 \ar @{.>} "A";"B" <0pt>
\endxy} \Ea
\in \RGra_d(1,2), \ \
\Ba{c}\resizebox{17mm}{!}{
\mbox{$\xy
 (-2,3)*{^{\bar{1}}};
 (-2,-3)*{_{\bar{2}}};
 (-9,0)*{^{\bar{3}}};
 (-4,0)*+{_{_4}}*\frm{o}="C";
  (9,0)*+{_{_1}}*\frm{o}="1";
(-7,8)*+{_{_2}}*\frm{o}="2";
(-7,-8)*+{_{_3}}*\frm{o}="3";
 \ar @{.>} "1";"C" <0pt>
  \ar @{.>} "1";"2" <0pt>
   \ar @{.>} "1";"3" <0pt>
 \ar @{.>} "2";"C" <0pt>
  \ar @{.>} "3";"C" <0pt>
\endxy$}}
\Ea
\in \RGra_d(3,4).
$$
Each edge of $\Ga$ is equipped with a direction which can be flipped producing the
following sign factor,
\Beq\label{2: symmetry of dotted edges}
\Ba{c} \resizebox{11mm}{!}{\xy
 (10,1)*+{_k}*\frm{o}="B";
 (0,1)*+{_i}*\frm{o}="A";
 \ar @{.>} "A";"B" <0pt>
\endxy} \Ea
= (-1)^{d}
\Ba{c} \resizebox{11mm}{!}{\xy
 (10,1)*+{_k}*\frm{o}="B";
 (0,1)*+{_i}*\frm{o}="A";
 \ar @{<.} "A";"B" <0pt>
\endxy} \Ea.
\Eeq
Hence one can skip from now on showing directions on dotted edges (assuming tacitly that some choice has been made). The cohomological degree of $\Ga$ is 
defined by
$$
|\Ga|=(1-d) \# E(\Ga).
$$
 For
for $d$ even it is tacitly  assumed that some ordering of edges of $\Ga$ is fixed up to an even 
permutation (an odd permutation acts as the multiplication by $-1$), while for $d$ odd it is 
assumed that some direction on each dotted edge is chosen (up to a flip as in 
(\ref{2: symmetry of dotted edges})).

\sip

The (reduced) properadic compositions in $\RGra_d$ 
$$
 \Ba{rccc}
 _i\circ_j:   & \cR\cG ra_d(m_1,n_1)\ot_\K \cR\cG ra_d(m_2,n_2) &\lon & 
 \cR\cG ra_d(m_1+m_2-1,n_1+n_2-1)\\
        &    \Ga_1\ot  \Ga_2 & \lon & {\Ga_1}\  _i\circ_j \Ga_2
        \Ea
 $$
are given, for any $i\in [n_1]$ and any $j\in [m_2]$, by substituting the 
 boundary $b_j$ of the ribbon graph $\Ga_2$ into the vertex $v_i$  of the
 ribbon graph $\Ga_1$, and then re-attaching  half-edges (glued earlier 
to $v_i$) to corners of the boundary $b_j$ in all possible ways while respecting 
the cyclic orders in the sets $H(v_i)$ and $C(b_j)$; put another way, we sum over all possible maps
$H(v_i)\rar C(b_j)$ of cyclic sets (see \S 4.2 in 
\cite{MW1} for full details and illustrating examples). The composition is best understood using 
pictorial representations of the boundary and of the vertex,
\Beq\label{2: composition in RGra}
b_j=  \Ba{c}\resizebox{15mm}{!}{\xy
(6,-8)*{\circlearrowleft},
 (0,0)*{}="a1",
(10,0)*{}="a2",
(13,-3)*{}="a3",
(13,-13)*{}="a4",
(10,-16)*{}="a5",
(0,-16)*{}="a6",
(-3,-13)*{}="a7",
(-3,-3)*{}="a8",
\ar @{-} "a1";"a2" <0pt>
\ar @{.} "a2";"a3" <0pt>
\ar @{-} "a3";"a4" <0pt>
\ar @{.} "a4";"a5" <0pt>
\ar @{-} "a5";"a6" <0pt>
\ar @{.} "a6";"a7" <0pt>
\ar @{-} "a7";"a8" <0pt>
\ar @{.} "a8";"a1" <0pt>
\endxy}\Ea, \ \ \ \ 
v_i= \Ba{c}\resizebox{21mm}{!}{\xy
 (0,0)*{
\xycircle(6,6){.}};
(-6,0)*{}="1";
(-16,0)*{}="1'";
(4.2,4.2)*{}="2";
(12,12)*{}="2'";
(4.2,-4.2)*{}="3";
(12,-12)*{}="3'";
(6,0)*{}="4";
(16,0)*{}="4'";
(-4.2,4.2)*{}="5";
(-12,12)*{}="5'";
(-4.2,-4.2)*{}="6";
(-12,-12)*{}="6'";
(-0.3,6)*{}="a";
(-0.4,6)*{}="b";
\ar @{-} "1";"1'" <0pt>
\ar @{-} "2";"2'" <0pt>
\ar @{-} "3";"3'" <0pt>
\ar @{-} "4";"4'" <0pt>
\ar @{-} "5";"5'" <0pt>
\ar @{->} "a";"b" <0pt>
\ar @{-} "6";"6'" <0pt>
\endxy}\Ea,
\ \ \ 
_i\circ_j: \ \sum
 \Ba{c}\resizebox{20mm}{!}{\xy
(6,-8)*{\circlearrowleft},
 (0,0)*{}="a1",
(10,0)*{}="a2",
(13,-3)*{}="a3",
(13,-13)*{}="a4",
(10,-16)*{}="a5",
(0,-16)*{}="a6",
(-3,-13)*{}="a7",
(-3,-3)*{}="a8",
(-1,-1)*{}="1";
(-6,7)*{}="1'";
(12,-2)*{}="2";
(20,4)*{}="2'";
(12,-14)*{}="3";
(18,-19)*{}="3'";
(-1.8,-1.8)*{}="4";
(-8,5)*{}="4'";
(-2.5,-2.5)*{}="5";
(-11,2)*{}="5'";
(11,-15)*{}="6";
(15,-22)*{}="6'";
\ar @{-} "1";"1'" <0pt>
\ar @{-} "2";"2'" <0pt>
\ar @{-} "3";"3'" <0pt>
\ar @{-} "4";"4'" <0pt>
\ar @{-} "5";"5'" <0pt>
\ar @{-} "6";"6'" <0pt>
\ar @{-} "a1";"a2" <0pt>
\ar @{.} "a2";"a3" <0pt>
\ar @{-} "a3";"a4" <0pt>
\ar @{.} "a4";"a5" <0pt>
\ar @{-} "a5";"a6" <0pt>
\ar @{.} "a6";"a7" <0pt>
\ar @{-} "a7";"a8" <0pt>
\ar @{.} "a8";"a1" <0pt>
\endxy}
\Ea
\Eeq

There is a  morphism of properads \cite{MW1}
\Beq\label{2: i from Lie to RGra}
i: \Lie_d \lon    \RGra_d
\Eeq
given on the Lie bracket generator of $\Lie_d$  by
$$
\Ba{c}\begin{xy}
 <0mm,0.66mm>*{};<0mm,3mm>*{}**@{-},
 <0.39mm,-0.39mm>*{};<2.2mm,-2.2mm>*{}**@{-},
 <-0.35mm,-0.35mm>*{};<-2.2mm,-2.2mm>*{}**@{-},
 <0mm,0mm>*{\bu};<0mm,4.1mm>*{^{^{\bar{1}}}}**@{},
   <0.39mm,-0.39mm>*{};<2.9mm,-4mm>*{^{_2}}**@{},
   <-0.35mm,-0.35mm>*{};<-2.8mm,-4mm>*{^{_1}}**@{},
\end{xy}\Ea
 \stackrel{i}{\lon}
 \Ba{c} \resizebox{11mm}{!}{\xy
 (10,1)*+{_2}*\frm{o}="B";
 (0,1)*+{_1}*\frm{o}="A";
 \ar @{.} "A";"B" <0pt>
\endxy} \Ea
$$
Applying to $\RGra_d$ Thomas Willwacher's  \cite{Wi1}
twisting endofunctor $\tw$  one obtains a dg properad
$$
\tw \RGra_d=\left\{ \tw \RGra_d(m,n), \delta   \right\}_{m\geq 1, n \geq 0}
$$
which is generated by ribbon graphs $\Ga$ with $m\geq 1$ labelled boundaries, $n\geq 0$ labelled vertices, and any number of unlabelled vertices $\circ$  to which one assigns 
the cohomological degree $d$ (see \S 3.3 in \cite{Me1} for a more detailed description). 
For example,
$$
\Ba{c} \resizebox{11mm}{!}{\xy
 (5,3)*{^{\bar{1}}};
 (10,1)*+{_1}*\frm{o}="B";
 (0,1)*+{\circ}="A";
 \ar @{.} "A";"B" <0pt>
\endxy} \Ea
\in \tw \RGra_d(1,1), \ \
\Ba{c}\resizebox{17mm}{!}{
\mbox{$\xy
 (-2,3)*{^{\bar{1}}};
 (-2,-3)*{_{\bar{2}}};
 (-9,0)*{^{\bar{3}}};
 (-4,0)*{\circ}="C";
  (9,0)*+{_{_1}}*\frm{o}="1";
(-7,8)*+{_{_2}}*\frm{o}="2";
(-7,-8)*{\circ}="3";
 \ar @{.} "1";"C" <0pt>
  \ar @{.} "1";"2" <0pt>
   \ar @{.} "1";"3" <0pt>
 \ar @{.} "2";"C" <0pt>
  \ar @{.} "3";"C" <0pt>
\endxy$}}
\Ea
\in \tw \RGra_d(3,2), \ \
\Ba{c}\resizebox{12mm}{!}{\xy
(0,0)*{_{\bar{1}}};
 (0,5)*{_{\bar{2}}};
   {\ar@{.}@/^0.7pc/(-6,0)*{\circ} ;(6,0)*{\circ}};
 {\ar@{.}@/^0.7pc/(6,0)*{\circ};(-6,0)*{\circ}};
\endxy}\Ea\in \tw \RGra_d(2,0).
$$
The cohomological degree of $\Ga\in \tw\RGra_d$ is given by
$$
|\Ga|=(1-d) \# E(\Ga) + d\# V_\circ(\Ga)
$$
where $V_\circ(\Ga)$ stands for the set of unlabelled vertices.
The differential in $\tw\RGra_d(m,n)$ is given by the standard ``splitting of vertices" formula,
\Beq\label{2: delta in Tw(RGra)}
\delta\Ga:=
\sum_{i=1}^m \Ba{c}\resizebox{4mm}{!}{  \xy
 (0,7)*{\circ}="A";
 (0,0)*+{_1}*\frm{o}="B";
 \ar @{.} "A";"B" <0pt>
\endxy} \Ea
\  _1\hspace{-0.7mm}\circ_i \Ga\ \
- \ \ (-1)^{|\Ga|} 
\sum_{j=1}^n\Ga\  _j\hspace{-0.7mm}\circ_1
\Ba{c}\resizebox{4mm}{!}{  \xy
 (0,6)*{\circ}="A";
 (0,0)*+{_1}*\frm{o}="B";
 \ar @{.} "A";"B" <0pt>
\endxy} \Ea
\  -(-1)^{|\Ga|}\ \frac{1}{2} \sum_{v\in V_\circ(\Ga)} \Ga\circ_v  \left(\xy
 (0,0)*{\circ}="a",
(5,0)*{\circ}="b",
\ar @{.} "a";"b" <0pt>
\endxy\right), \Eeq
where the symbol $_i\circ_j$ stands for the properadic composition as in (\ref{2: composition in RGra}), and  the symbol $\Ga\circ_v  \left(\xy
 (0,0)*{\circ}="a",
(5,0)*{\circ}="b",
\ar @{.} "a";"b" <0pt>
\endxy\right)$ means the substitution of the graph  $\xy
 (0,0)*{\circ}="a",
(5,0)*{\circ}="b",
\ar @{.} "a";"b" <0pt>
\endxy$ into the unlabelled vertex $v$ of the graph $\Ga$ followed by
 the summation over all possible re-attachments of the half-edges attached earlier
 to $v$ among the two unlabelled vertices in all possible ways which respect their cyclic
 orderings (again in a full analogy to  (\ref{2: composition in RGra})). Note that for almost all graphs the univalent unlabelled vertex
 created in the first summand of $\delta$ cancel out the univalent
 unlabelled vertices created in the second and the third parts of that differential. If $\Ga$ has at least one edge and
 has no univalent unlabelled vertices, then $\delta\Ga$ will not have univalent vertices
 either.

\sip

 The dg properad  $\tw\cR\cG ra_d$ contains a dg sub-properad $\RGraphs_d$
  spanned by ribbon graphs having at least one unlabelled vertex of valency  $\geq 3$ (it is not a big loss to work with $\RGraphs_d$ instead of the full properad $\tw\cR\cG ra_d$, see \S 3.4 of
  \cite{Me1} for the explicit interrelation between the associated cohomology groups).

\sip

As a complex, each $\bS_m^{op}\times \bS_n$-module $\RGraphs_d(m,n)$ decomposes into a direct
product,
$$
\RGraphs_d(m,n)=\prod_{g\geq 0 \atop 2g+m+n\geq 3} \RGraphs_d^{(g)}(m,n)
$$
where $\RGraphs_d^{(g)}(m,n)$ is generated by ribbon graphs $\Ga$ of genus $g$. One has for any $m\geq 1$, $n\geq 0$ \cite{Co,Me1}
$$
H^\bu(\RGraphs_d^{(g)}(m,n), \delta) =H_c^{\bu -m +d(2g-2+m+n) }(\cM_{g, m+n}), $$
 where $\cM_{g,m+n}$ is the moduli space of algebraic curves of genus $g$ with $m\geq 1$ 
 boundaries and   $n\geq 0$ marked points.
The complex 
\Beq\label{2: RGCdm}
 \RGC_d^{(m)}:= \prod_{g\geq 0} \RGC_d^{(g,m)}, \ \ \text{where}  \ \     \RGC_d^{(g,m)}:= \left(\RGraphs_d^{(g)}(m,0)[d], \delta\right),  \ \ \ 
\Eeq
 is (up to a degree shift) the classical R.\ Penner's ribbon graph complex \cite{Pe} 
 with marked boundaries  (see also \cite{Ko1}) which computes the compactly supported 
 cohomology,
$$
H^\bu(\RGC_d^{(g,m)}) = H_c^{\bu+ (d-1)m +d(2g-1)}(\cM_{g,m}),
$$
of the moduli spaces $\cM_{g,m}$  of genus $g$ algebraic curves with $m$ marked points. 

\subsection{A dg Lie algebra of ribbon graphs with one boundary}\label{2: subsec on rgc} Let
$$
ope\RGra_d:=\{ \RGra_d(1,n)\}_{n\geq 1} \subset \RGra_d
$$
be a suboperad of the properad $\RGra_d$ generated by ribbon graphs with precisely one boundary.
The morphism (\ref{2: i from Lie to RGra}) factors through the morphism of operads
$$
i: \Lie_d \lon    ope\RGra_d.
$$
Hence one can consider the deformation complex of $i$ in the category of operads \cite{L} 
$$
\rgc_d:= \Def(\Lie_d \stackrel{i}{\rar}ope\RGra_d)\simeq \prod_{n\geq 1} 
\RGra_d(1,n)\ot_{\bS_n} \sgn_n^{|d|}[d-dn],
$$
which has a canonical pre-Lie algebra structure given by
$$
\Ba{rccc}
\circ: & \rgc_d\ot \rgc_d & \lon & \rgc_d\\
       &  \Ga_1 \ot \Ga_2 & \lon & \Ga_1\circ \Ga_2:= 
       \sum_{v\in V(\Ga_1)} {\Ga_1}\ _v\hspace{-0.8mm}\circ_b \Ga_2
\Ea
$$
where $_v\circ_b$ stand for the standard substitution of the unique boundary $b$ of $\Ga_2$
into the vertex $v$ of $\Ga_1$ as in (\ref{2: composition in RGra}). The Lie bracket (of degree zero) in
$\rgc_d$ is given by
$$
[\Ga_1,\Ga_2]= \Ga_1\circ \Ga_2 - (-1)^{|\Ga_1||\Ga_2|} \Ga_2\circ\Ga_1.
$$
This (pre)Lie algebra is graded with respect to the genus $g$ of the generating 
ribbon graphs. 
The differential in the deformation complex $\rgc_d$ is given by
$$
\delta\Ga:= \left[ \xy
 (0,0)*{\circ}="a",
(5,0)*{\circ}="b",
\ar @{.} "a";"b" <0pt>
\endxy, \Ga\right].
$$
The cohomological degree of $\Ga\in \rgc_d$ is given by
$$
|\Ga|=d(\#V(\Ga)-1)+(1-d)\# E(\Ga).
$$
There is an obvious isomorphism of complexes (which explains the degree shift in the definition
of $\RGC_d^{g,m}$ above)
$$
\rgc_d= \prod_{g\geq 1} \RGC_d^{g,1}
$$
so that the cohomology of the dg Lie algebra  $\rgc_d$ is given by
\Beq\label{2: H(rgc_d)}
H^\bu(\rgc_d)=\prod_{g\geq 1}  H_c^{\bu+2gd-1}(\cM_{g,1}).
\Eeq

\subsection{A dg properad  of ribbon quivers}\label{2: subsec on RGra^or and RGraphs^or}
Let $\DRGra_{d}=\{\DRGra_{d}(m,n)\}_{m,n\geq 1}$ be a version of  the properad $\RGra_{d}$ 
which is generated by ribbon graphs with the condition (\ref{2: symmetry of dotted edges}) on edges 
dropped, that is, every  edge comes equipped with a fixed direction. In our pictures we show such edges  as {\em solid arrows}, e.g.
 $$
\Ba{c} \resizebox{11mm}{!}{\xy
 (5,3)*{^{\bar{1}}};
 (10,1)*+{_1}*\frm{o}="B";
 (0,1)*+{_2}*\frm{o}="A";
 \ar @{->} "A";"B" <0pt>
\endxy} \Ea
\in \DRGra_d(1,2), \ \
\Ba{c}\resizebox{17mm}{!}{
\mbox{$\xy
 (-5,2)*{^{\bar{1}}};
 (-5,-3)*{_{\bar{2}}};
 (-12,4)*{^{\bar{3}}};
 (-12,0)*+{_{_4}}*\frm{o}="C";
  (4,0)*+{_{_1}}*\frm{o}="1";
(-7,8)*+{_{_2}}*\frm{o}="2";
(-7,-8)*+{_{_3}}*\frm{o}="3";
 \ar @{->} "1";"C" <0pt>
  \ar @{->} "1";"2" <0pt>
   \ar @{->} "1";"3" <0pt>
 \ar @{->} "2";"C" <0pt>
  \ar @{>} "3";"C" <0pt>
\endxy$}}
\Ea
\in \DRGra_d(3,4)
, \ \
\Ba{c}\resizebox{17mm}{!}{
\mbox{$\xy
 (-5,2)*{^{\bar{1}}};
 (-5,-3)*{_{\bar{2}}};
 (-12,4)*{^{\bar{3}}};
 (-12,0)*+{_{_4}}*\frm{o}="C";
  (4,0)*+{_{_1}}*\frm{o}="1";
(-7,8)*+{_{_2}}*\frm{o}="2";
(-7,-8)*+{_{_3}}*\frm{o}="3";
 \ar @{<-} "1";"C" <0pt>
  \ar @{->} "1";"2" <0pt>
   \ar @{->} "1";"3" <0pt>
 \ar @{->} "2";"C" <0pt>
  \ar @{>} "3";"C" <0pt>
\endxy$}}
\Ea
\in \DRGra_d(3,4).
 $$
This properad contains a sub-properad
$\ORGra_d$ generated by ribbon graphs with no closed paths of directed edges;
for example, all the above ribbon graphs except the last one belong to $\ORGra_{d}$. Such graphs
 are often called {\em oriented}\, (cf.\ \cite{Wi2}); in this paper we also call them  {\em ribbon quivers}.

\sip

There is a morphism of properads
\Beq\label{2: map from lie to RGra_or}
i^{or}: \Lie_d \lon   \ORGra_d
\Eeq
given on the Lie bracket generator of $\Lie_d$ by
$$
i^{or}:
\Ba{c}\begin{xy}
 <0mm,0.66mm>*{};<0mm,3mm>*{}**@{-},
 <0.39mm,-0.39mm>*{};<2.2mm,-2.2mm>*{}**@{-},
 <-0.35mm,-0.35mm>*{};<-2.2mm,-2.2mm>*{}**@{-},
 <0mm,0mm>*{\bu};<0mm,4.1mm>*{^{^{\bar{1}}}}**@{},
   <0.39mm,-0.39mm>*{};<2.9mm,-4mm>*{^{_2}}**@{},
   <-0.35mm,-0.35mm>*{};<-2.8mm,-4mm>*{^{_1}}**@{},
\end{xy}\Ea
 \lon
 \frac{1}{2}\left(
 \Ba{c}\resizebox{11mm}{!}{  \xy
 (3.5,4)*{^{\bar{1}}};
 (9,1)*+{_2}*\frm{o}="A";
 (0,1)*+{_1}*\frm{o}="B";
 \ar @{<-} "A";"B" <0pt>
\endxy} \Ea
 +(-1)^d
 \Ba{c}\resizebox{11mm}{!}{  \xy
 (3.5,4)*{^{\bar{1}}};
 (9,1)*+{_1}*\frm{o}="A";
 (0,1)*+{_2}*\frm{o}="B";
 \ar @{<-} "A";"B" <0pt>
\endxy} \Ea \right).
$$
Applying the twisting endofunctor \cite{Wi1} to $\ORGra_d$ one obtains, in a full analogy to $\tw\RGra_d$ above,  a {\it differential}\,  
graded properad (cf.\ \cite{Me3})
$$
\tw \ORGra_d=:\left\{\tw \ORGra_d(m,n), \delta   \right\}_{m\geq 1, n \geq 0}
$$
which is generated by ribbon quivers $\Ga$ with $m\geq 1$ labelled boundaries, $n\geq 0$ 
labelled vertices, and any
number of unlabelled vertices which are shown in pictures in black color as $\bu$,
 and which are assigned the cohomological degree $d$. For example,
$$
\Ba{c} \resizebox{11mm}{!}{\xy
 (5,3)*{^{\bar{1}}};
 (10,1)*+{_1}*\frm{o}="B";
 (0,1)*{\bu}="A";
 \ar @{->} "A";"B" <0pt>
\endxy} \Ea
\in \tw \ORGra_d(1,1), \ \
\Ba{c}\resizebox{17mm}{!}{
\mbox{$\xy
 (-2,3)*{^{\bar{1}}};
 (-2,-3)*{_{\bar{2}}};
 (-9,0)*{^{\bar{3}}};
 (-4,0)*{\bu}="C";
  (9,0)*+{_{_1}}*\frm{o}="1";
(-7,8)*+{_{_2}}*\frm{o}="2";
(-7,-8)*{\bu}="3";
 \ar @{->} "1";"C" <0pt>
  \ar @{->} "1";"2" <0pt>
   \ar @{<-} "1";"3" <0pt>
 \ar @{->} "2";"C" <0pt>
  \ar @{->} "3";"C" <0pt>
\endxy$}}
\Ea
\in\tw \ORGra_d(3,2), \ \
\Ba{c}\resizebox{12mm}{!}{\xy
(0,0)*{_{\bar{1}}};
 (0,5)*{_{\bar{2}}};
   {\ar@{->}@/^0.7pc/(-6,0)*{\bu} ;(6,0)*{\bu}};
 {\ar@{<-}@/^0.7pc/(6,0)*{\bu};(-6,0)*{\bu}};
\endxy}\Ea\in \tw\ORGra_d(2,0).
$$
The cohomological degree of $\Ga\in \ORGraphs_d$ is given by the formula
$$
|\Ga|=(1-d) \# E(\Ga) + d\# V_\bu(\Ga)
$$
where $V_\bu(\Ga)$ stands for the set of unlabelled vertices.
The differential in $\tw\ORGra_d(m,n)$ is given by the standard ``splitting of vertices" formula (cf.\ (\ref{2: delta in Tw(RGra)}))
\Beq\label{2: delta in Tw(RGra^or)}
\delta\Ga:=
\sum_{i=1}^m \left(\hspace{-2mm} \Ba{c}\resizebox{4mm}{!}{  \xy
 (0,7)*{\bu}="A";
 (0,0)*+{_1}*\frm{o}="B";
 \ar @{->} "A";"B" <0pt>
\endxy} \Ea
+(-1)^d\hspace{-2mm}
\Ba{c}\resizebox{4mm}{!}{  \xy
 (0,7)*{\bu}="A";
 (0,0)*+{_1}*\frm{o}="B";
 \ar @{<-} "A";"B" <0pt>
\endxy} \Ea
\right)
\  _1\hspace{-0.7mm}\circ_i \Ga\ \
- \ \ (-1)^{|\Ga|} 
\sum_{j=1}^n\Ga\  _j\hspace{-0.7mm}\circ_1
\left(\hspace{-2mm} \Ba{c}\resizebox{4mm}{!}{  \xy
 (0,7)*{\bu}="A";
 (0,0)*+{_1}*\frm{o}="B";
 \ar @{->} "A";"B" <0pt>
\endxy} \Ea
+(-1)^d\hspace{-2mm}
\Ba{c}\resizebox{4mm}{!}{  \xy
 (0,7)*{\bu}="A";
 (0,0)*+{_1}*\frm{o}="B";
 \ar @{<-} "A";"B" <0pt>
\endxy} \Ea
\right)
\  -\hspace{-2mm}\sum_{v\in V_\bu(\Ga)} \Ga\circ_v  \left(\xy
 (0,0)*{\bu}="a",
(5,0)*{\bu}="b",
\ar @{->} "a";"b" <0pt>
\endxy\right).
\Eeq
Every generator $\Ga$ of $\tw\ORGra(m,n)$ comes equipped with an orientation 
$or(\Ga)$ which is, by definition, a unital vector{\footnote{For a finite set $S$ let us denote the top degree skew-symmetric tensor power
of $\K\left\langle S\right\rangle$ by $\det S$. Let us assume that $\det S$
is a 1-dimensional Euclidean space associated with the unique Euclidean
structure on $\K\left\langle S\right\rangle$ in which the elements of $S$ serve
as an orthonormal basis; in particular, $\det S$ contains precisely two vectors of unit length.
Then $or(\Ga)$ is, by definition a unital vector in $\det V(\Ga)$ (resp., in $\det E(\Ga)$) for $d$ odd 
(resp., for $d$ even).} in $\det V_\bu(\Ga)$ (resp., in $\det E(\Ga)$) for $d$ odd (resp., $d$ even).

\sip

Given any pre-CY structure $\pi$ in $A$,  there is an associated action 
of the dg properad $\tw \RGra^{or}_d$ on the extended higher Hochschild complex
$$
\widehat{C}_{[d]}^\bu(A) :=\prod_{k\geq 0} C_{[d]}^{(k)}(A):= \prod_{k\geq 0}
\left( \bigoplus_{n_1,\ldots, n_k\geq 0}
\Hom\left(\bigotimes_{i=1}^{k}
(A[1])^{\ot n_i}, (A[2-d])^{\ot k}\right)_{\Z_k}\right)
$$
of $(A,\pi)$ (see \cite{Me3} for more details). This properad was used in \cite{Me3} to 
construct a so called 
oriented gravity properad $\OORGraphs_d$ which has some useful properties in the context of 
the theory of pre-CY algebras. 

\sip

We define a subcomplex
 \Beq\label{2: ORGC_d,m}
 \ORGC_d^{(m)} \subset \tw\RGra_d^{or}(m,0)[d]
 \Eeq
generated by ribbon graphs having at least 
one black vertex of valency $\geq 3$. The above inclusion is a quasi-isomorphism up to the standard
series of ribbon quivers with unlabelled {\em bivalent}\, vertices (cf.\ \cite{Wi1,Wi2,Me1}).


\subsection{A dg Lie algebra of ribbon quivers with one boundary}\label{2: subsec on orgc} This subsection is a straightforward quiver version of the story presented in \S {\ref{2: subsec on rgc}}. Let
$$
ope\ORGra_d:=\{ \ORGra_d(1,n)\}_{n\geq 1} \subset \ORGra_d
$$
be a suboperad of the properad $\ORGra_d$ generated by ribbon quivers with precisely one boundary.
The morphism (\ref{2: map from lie to RGra_or}) factors through the morphism of {\em operads}\,
$$
i^{or}: \Lie_d \lon    ope\ORGra_d.
$$
Hence one can consider the deformation complex of $i^{or}$ in the category of operads \cite{L,MV}
$$
\orgc_d:= \Def(\Lie_d \stackrel{i}{\rar}ope\ORGra_d)\simeq \prod_{n\geq 1} \ORGra_d(1,n)\ot_{\bS_n} \sgn_n^{|d|}[d-dn]
$$
which has a canonical pre-Lie algebra structure $\circ$ which is given by substitunig the unique boundary 
of one graph into vertices of another graph, in a full analogy to the case of $\rgc_d$.
The differential in the pre-Lie algebra $\orgc_d$ is given by
$$
\delta\Ga:= \left[ \xy
 (0,0)*{\bu}="a",
(5,0)*{\bu}="b",
\ar @{->} "a";"b" <0pt>
\endxy, \Ga\right].
$$
The ribbon graphs generating  $\orgc_d$ must have at least one vertex of valency $\geq 3$.

\mip

{\Large
\section{\bf Proof of Theorem {\ref{1: Main theorem 1}}}
}

\sip

\subsection{An auxiliary dg properad $\OORGraphs_{d,d+1}$}
Let us consider again the twisted properad $\tw\ORGra_{d+1}$ from \S {\ref{2: subsec on RGra^or and RGraphs^or}} but equipped now not with the full differential $\delta$ as in  (\ref{2: delta in Tw(RGra^or)}), but with its much abbreviated version 
$$
\delta_\bu (\Ga):= -\hspace{-2mm}\sum_{v\in V_\bu(\Ga)} \Ga\circ_v  \left(\xy
 (0,0)*{\bu}="a",
(5,0)*{\bu}="b",
\ar @{->} "a";"b" <0pt>
\endxy\right).
$$
The point is that the resulting dg properad
$(\tw\ORGra_{d+1},
\delta_\bu)$ contains a {\em dg}\, sub-properad $(\widehat{O}\RGraphs_{d+1}, \delta_\bu)$ which is, by definition, generated by ribbon
quivers with no outgoing edges at labelled  vertices (this not true for $(\tw\ORGra_{d+1},\delta)$ as the full differential $\delta$ does not preserve this condition).
It was noticed in \cite{Me3} that there is an explicit morphism
\Beq\label{3: morphism f from cHoile}
f: (c\Holie_d, \delta) \lon  (\widehat{O}\RGraphs_{d+1},
\delta_\bu)
\Eeq
from the dg operad $c\Holie_d$  of curved strongly homotopy Lie algebras 
which is given on the generators by 
$$
\Ba{c}\resizebox{22mm}{!}{ \xy
(1,-5)*{\ldots},
(-13,-7)*{_1},
(-8,-7)*{_2},
(-3,-7)*{_3},
(7,-7)*{_{n-1}},
(13,-7)*{_n},
 (0,0)*{\bu}="a",
(0,5)*{}="0",
(-12,-5)*{}="b_1",
(-8,-5)*{}="b_2",
(-3,-5)*{}="b_3",
(8,-5)*{}="b_4",
(12,-5)*{}="b_5",
\ar @{-} "a";"0" <0pt>
\ar @{-} "a";"b_2" <0pt>
\ar @{-} "a";"b_3" <0pt>
\ar @{-} "a";"b_1" <0pt>
\ar @{-} "a";"b_4" <0pt>
\ar @{-} "a";"b_5" <0pt>
\endxy}\Ea
\lon
\sum_{\sigma\in \bS_n\atop i_k:=\sigma(k)}\frac{(-1)^{d\sigma}}{n}\
{\resizebox{21mm}{!}{
\xy
(-6,-9)*+{_{_{i_n}}}*\frm{o}="a1",
(-6,9)*+{_{_{i_2}}}*\frm{o}="a2",
(6,9)*+{_{_{i_3}}}*\frm{o}="a3",
(10,0)*+{_{_{i_4}}}*\frm{o}="a4",
(6,-9)*+{_{...}}*\frm{o}="a5",
(-10,0)*+{_{_{i_1}}}*\frm{o}="a6";
(0,0)*{\bu}="b",
\ar @{<-} "a1";"b" <0pt>
\ar @{<-} "a2";"b" <0pt>
\ar @{<-} "a3";"b" <0pt>
\ar @{<-} "a4";"b" <0pt>
\ar @{<-} "a5";"b" <0pt>
\ar @{<-} "a6";"b" <0pt>
\endxy}} \ \ \  \ \ \  \forall  n\geq 0.
$$
Applying the standard twisting endofunctor \cite{Wi1, DSV} to the above map  one obtains a  dg properad   $\tw\widehat{O}\RGraphs_{d+1}$  which played a central role in the proofs of main statements in \cite{Me3}. It is
generated by ribbon quivers $\Ga$ with two types of unlabelled vertices, black ones of degree $d+1$
and white ones of degree $d$, 
$$
\Ba{c}\resizebox{13mm}{!}{ \xy
(-2.5,2.5)*{^{\bar{1}}},
(-7,4)*{^{\bar{2}}},
(3,5)*{^{\bar{3}}},
(-5,0)*+{_1}*\frm{o}="1";
 (5,0)*{\bu}="11";
 (5,10)*{\circ}="2";
 (-5,10)*{\bu}="22";
 \ar @{<-} "1";"22" <0pt>
 \ar @{->} "11";"22" <0pt>
 \ar @{->} "11";"22" <0pt>
  \ar @{->} "11";"1" <0pt>
   \ar @{->} "11";"2" <0pt>
  \ar @{<-} "2";"22" <0pt>
\endxy} \Ea \in  \tw\widehat{O}\RGraphs_{d+1}(3,1),
$$
such that both labelled and unlabelled {\it white}\, vertices have no outgoing edges. The cohomological degree of such a quiver is given by
$$
|\Ga|=(d+1)\# V_\bu(\Ga) + d\# V_\circ(\Ga) - d \# E(\Ga),
$$
while the differential in  $\tw\widehat{O}\RGraphs_{d+1}$ is given by the sum \cite{Me3}
$$
\delta \Ga:= d_\bu \Ga + D_\ga(\Ga).
$$
where 
$$
d_\bu \Ga:=- \left(\sum_{v\in V_{\circ}(\Ga)}\delta_{\circ\bu}^v\Ga + \Ga \circ_v \left( \sum_{k=1}^\infty \frac{1}{k}
{\resizebox{24mm}{!}{
\xy
(-4,-5)*{\circ}="a1",
(-4,5)*{\circ}="a2",
(4,5)*{\circ}="a3",
(4,-5)*{\circ}="a4",
(7,0)*{\circ}="a5",
(-7,0)*{\circ}="a6",
(0,0)*{\bu}="b",
(6,-9)*{\underbrace{\hspace{12mm}}_{k\ \text{edges}} \hspace{12mm}},
\ar @{<-} "a1";"b" <0pt>
\ar @{<-} "a2";"b" <0pt>
\ar @{<-} "a3";"b" <0pt>
\ar @{<-} "a4";"b" <0pt>
\ar @{<-} "a5";"b" <0pt>
\ar @{<-} "a6";"b" <0pt>
\endxy}}\right)
+
\sum_{w\in V_{\bu}(\Ga) } \Ga \circ_w \left( \resizebox{2.6mm}{!}{
\xy
 (0,4)*{\bullet}="a",
(0,-2)*{\bu}="b",
\ar @{->} "a";"b" <0pt>
\endxy}
\right)\right)
$$
and 
$$
D_\ga(\Ga):=\sum_{i=1}^m \ga\, {_1\circ_i}\Ga - (-1)^{|\Ga||\ga|} \sum_{j=1}^n \Ga{_j\circ_1}\ga.
\ \ \ \ga:= \sum_{k=0}^\infty
{\resizebox{24mm}{!}{
\xy
(-4,-5)*{\circ}="a1",
(-4,5)*{\circ}="a2",
(4,5)*{\circ}="a3",
(4,-5)*{\circ}="a4",
(7,0)*{\circ}="a5",
(-7,0)*+{_1}*\frm{o}="a6";
(0,0)*{\bu}="b",
(6,-9)*{\underbrace{\hspace{12mm}}_{k+1\ \text{edges}} \hspace{12mm}},
\ar @{<-} "a1";"b" <0pt>
\ar @{<-} "a2";"b" <0pt>
\ar @{<-} "a3";"b" <0pt>
\ar @{<-} "a4";"b" <0pt>
\ar @{<-} "a5";"b" <0pt>
\ar @{<-} "a6";"b" <0pt>
\endxy}}
$$
Here $\delta_{\circ\bu}^v \Ga$ means the ribbon graph obtained from $\Ga$ by making the white 
unlabelled vertex $v$ into a black unlabelled vertex, and the symbol $\circ_v (X)$ (resp.\ $\circ_w (X)$)
 means the substitution of the unique boundary of the ribbon graph $X$ into the unlabelled white vertex $v$ (resp., black vertex $w$) and performing the properadic-like composition as in (\ref{2: composition in RGra}).

\sip

 We denote by
$$
 \OORGraphs_{d,d+1}=\left\{(\OORGraphs_{d,d+1}(m,n),\delta)\right\}_{m\geq 1,n\geq 0}.
$$
a dg sub-properad of $\tw\widehat{O}\RGraphs_{d+1}$ generated by ribbon graphs   having  
at least one unlabelled vertex (of any type) with valency $\geq 3$ and 
having no univalent unlabelled vertices.

\sip

Let us call a black bivalent vertex $\bu$ {\it special}\, if it has two outgoing edges 
attached, i.e.\ if it is of the form
\Beq\label{3: inessential vertex}
\xy
 (0,0)*{}="a",
(4,3)*{\bullet}="b",
(8,0)*{}="c",
\ar @{<-} "a";"b" <0pt>
\ar @{->} "b";"c" <0pt>
\endxy \ .
\Eeq

It was shown in \cite{Me3} (see  Proposition 4.4 there) that there is an explicit quasi-isomorphism of complexes
\Beq\label{3: pi from Rgraphs d,d+1 to RGrapgsd}
\pi_{m,n}: \OORGraphs_{d,d+1}(m,n) \lon \RGraphs_d(m,n), \ \ \ \ \forall\, m,n \geq 1
\Eeq
given by setting to zero every ribbon quiver in $\OORGraphs_{d,d+1}$ which has at least one 
{\it non-special black vertex}. 
Generators of the quotient dg properad get identified with generators
of $\RGraphs_d$ via the following trick,
\Beq\label{3: from special bu to dotted edges}
\Ba{c}\resizebox{18mm}{!}{
\xy
%
 (0,0)*+{u}="a",
(8,0)*{\bullet}="b",
(16,0)*+{v}="c",
\ar @{<-} "a";"b" <0pt>
\ar @{->} "b";"c" <0pt>
\endxy}\Ea
\  \to \
\Ba{c}\resizebox{11mm}{!}{
\xy
%
 (0,0)*{u}="a",
(8,0)*{v}="b",
\ar @{.} "a";"b" <0pt>
\endxy}\Ea
\Eeq
where $u$ and $v$ are arbitrary white vertices (labelled or unlabelled ones); 
the case $u=v$ is not excluded.
The argument in \cite{Me3} used heavily the fact that $n\geq 1$, i.e.\ that 
the generating graphs have at least one labelled vertex.
The situation with the case $n=0$ is quite different as the following Lemma shows.

\subsubsection{\bf Lemma}\label{3: Lemma on acyclicity of twORGraphs(0,m)}  {\it  The complex $\left(\OORGraphs_{d,d+1}(m,0), \delta\right)$
 is acyclic for any $m\geq 1$.}

\begin{proof} Consider a filtration of  $\left(\OORGraphs_{d,d+1}(m,0), \delta\right)$ by the total number of vertices. The induced differential $\delta_0$ in the associated graded $gr  \OORGraphs_{d,d+1}(m,0)$ acts on unlabelled white vertices
of ribbon quivers by making them black,
$$
\delta_0\Ga= -\sum_{v\in V_{\circ}(\Ga)}\delta_{\circ\bu}^v\Ga.
$$
To prove the acyclicity of $(gr  \OORGraphs_{d,d+1}(m,0), \delta_0)$ (and hence prove 
the Lemma), it is enough to prove the acyclicity of a version $gr  \OORGraphs^{marked}_{d,d+1}(m,0)$ of this complex in which all edges and vertices of the generators  are distinguished, but the type of previously unlabelled vertices is not fixed. If $G$ stands for the set of such generators, then 
there is an isomorphism of complexes  (cf.\ Lemma 6.2.1 in \cite{Me2}).
$$
gr  \OORGraphs^{marked}_{d,d+1}(m,0)=\prod_{\Ga\in G}\left(\bigotimes_{v\in V(\Ga)} C_v\right)
$$
where $C_v$ is a
\Bi
\item[(i)] one-dimensional trivial complex generated by one black vertex if $v$ has at least one outgoing solid edge in ${\Ga}$,
\item[(ii)] two-dimensional acyclic complex (generated by one black and one white vertex) if $v$ has no outgoing solid edges in $\Ga$.
\Ei
 As any generating
graph $\Ga$ has no closed paths of directed edges, then each $\Ga$  has at least one vertex of type (ii). Hence the complex 
$gr \OORGraphs^{marked}_{d,d+1}(m,0)$
has at least one acyclic tensor factor so that it is acyclic itself. The Lemma is proven.
\end{proof}

Let  $\OORGraphs_{d,d+1}^\bu (m,0)\hook \OORGraphs_{d,d+1}(m,0)$ be a subcomplex generated by ribbon quivers having at least one non-special black vertex. It fits a short exact
sequence of complexes
$$
0\lon \OORGraphs_{d,d+1}^\bu (m,0)[d] \lon \OORGraphs_{d,d+1}(m,0)[d] \lon \RGC_d^{(m)} \lon 0
$$
where $\RGC_d^{(m)}$ is R.\  Penner's ribbon graph complex (see \S 1.2).
The above Lemma implies

\subsubsection{\bf Corollary} 
$$
H^\bu\left(\RGC_d^{(m)}\right)= H^\bu\left( \OORGraphs_{d,d+1}^\bu (m,0)[d+1]\right).
$$

\sip

On the other hand, the complex  $\OORGraphs_{d,d+1}^\bu (m,0)[d+1]$ contains 
$\ORGC_{d+1}^{(m)}$ (defined in (\ref{2: ORGC_d,m})) as a subcomplex generated by 
ribbon quivers with no unlabelled white vertices.

\subsubsection{\bf Lemma}\label{3: Lemma on H(ORGCd) and d,d+1}
$$
H^\bu\left(\ORGC_{d+1}^{(m)}\right)= H^\bu\left( \OORGraphs_{d,d+1}^\bu (m,0)[d+1]\right).
$$
\begin{proof}
Consider a quotient complex $Z_{d,d+1}^{(m)}$ defined by the short exact sequence
\Beq\label{3: shoert exact seq with Zm}
0\lon \ORGC_{d+1}^{(m)} \lon \OORGraphs^\bu_{d,d+1}(m,0)[d+1] \lon Z_{d,d+1}^{(m)} \lon 0.
\Eeq
It is generated by ribbon quivers having at least one non-special black vertex  and at 
least one white vertex. To prove Lemma it is enough to show that $Z_{d,d+1}^{(m)}$ is acyclic.

\sip

Call a special black vertex $\bu$ {\em very special}\, if its two outgoing edges 
connect $\bu$ to some non-special black vertex $\bbu$ and some
white vertex $\circ$ as in the picture,
$$
\xy
%
 (0,0)*{\bbu}="a",
(8,0)*{\bullet}="b",
(16,0)*{\circ}="c",
\ar @{<-} "a";"b" <0pt>
\ar @{->} "b";"c" <0pt>
\endxy    .
$$
Consider a filtration of the complex $Z_{d,d+1}^{(m)}$ by the number black vertices 
which are not very special, and let $gr Z_{d,d+1}^{(m)}$ be the associated graded. 
To prove the acyclicity of $gr Z_{d,d+1}^{(m)}$ it is 
enough to prove the acyclicity of a version $\widehat{gr}Z_{d,d+1}^{(m)}$ of 
$gr Z_{d,d+1}$ in which all (black and white) vertices which are not very special are 
distinguished, say labelled by natural numbers.
Every generator $\Ga$ of $\widehat{gr}Z_{d,d+1}^{(m)}$
has at least non-special black vertex $\bbu$ (labelled say by $i$)  and at least one white vertex $\circ$ (labelled say by $k\neq i$) which are connected to each other either directly by an edge or via some very special vertex $\bu$ as in the following pictures  
$$
 \Ba{c}\resizebox{10mm}{!}{\xy
(0,2)*{^{i}},
(8,2)*{^{k}},
 (0,0)*{\bbu}="a",
 (8,0)*{\circ}="b";
\ar @{->} "a";"b" <0pt>
\endxy}\Ea, \ \ \ 
 \Ba{c}\resizebox{18mm}{!}{
\xy
(0,2)*{^{i}},
(16,2)*{^{k}},
 (0,0)*{\bbu}="a",
(8,0)*{\bullet}="b",
 (16,0)*{\circ}="c"
\ar @{<-} "a";"b" <0pt>
\ar @{->} "b";"c" <0pt>
\endxy}\Ea.
$$
We can assume without loss of generality that labels $k$ and $i$ take smallest possible values, say $i=1$ and $k=2$, and then consider a filtration of $\widehat{gr}Z_{d,d+1}^{(m)}$
by the number of very special vertices  which are connected to a non-special vertex with
label $\geq 3$. The associated graded
$gr \widehat{gr}Z_{d,d+1}^{(m)}$ is the tensor product of a trivial complex
and the complex $C_{12}$ which controls the types of all possible ``edges" between
vertices 1 and 2 (cf.\ Proposition 4.2 in \cite{MWW}). One has
$$
C_{12}=\bigoplus_{k\geq 1} \odot^k C
$$
where $C$ is a 2-dimensional complex generated by the following vectors,
$$
C=\text{span}\left\langle
 \Ba{c}\resizebox{10mm}{!}{\xy
(0,2)*{^{i}},
(8,2)*{^{k}},
 (0,0)*{\bbu}="a",
 (8,0)*{\circ}="b";
\ar @{->} "a";"b" <0pt>
\endxy}\Ea,
 \ \ \ 
 \Ba{c}\resizebox{18mm}{!}{
\xy
(0,2)*{^{i}},
(16,2)*{^{k}},
 (0,0)*{\bbu}="a",
(8,0)*{\bullet}="b",
 (16,0)*{\circ}="c"
\ar @{<-} "a";"b" <0pt>
\ar @{->} "b";"c" <0pt>
\endxy}\Ea
 \right\rangle
$$
and equipped with the differential sending the first vector to the second one.
This complex is acyclic implying the acyclicity of ${gr}Z_{d,d+1}^{(m)}$ 
and hence of $Z_{d,d+1}^{(m)}$.
\end{proof}

The above Corollary and Lemma imply the equality $H^\bu\left(\RGC_{d}^{(m)}\right)=H^\bu\left(\ORGC_{d+1}^{(m)}\right)$ for any $m\geq 1$.
Theorem {\ref{1: Main theorem 1}} is proven.

\mip

{\Large
\section{\bf Proof of Theorem {\ref{1: Main theorem 2}}}
}

\sip

\subsection{An auxiliary dg Lie algebra}
Let us consider the deformation complex of  the morphism 
(\ref{3: morphism f from cHoile}),
\Beqrn
\orgc_{d,d+1}&=& \Def\left(c\Holie_{d} \stackrel{f}{\rar}\widehat{O}\RGraphs_{d+1}\right)
   \simeq \prod_{n\geq 0} \widehat{O}\RGraphs_{d+1}(1,n)\ot_{\bS_n} \sgn_n^{|d|}[d-dn]\\
     &=&
    \underbrace{\widehat{O}\RGraphs_{d+1}(1,0)[d]}_{\simeq \orgc_{d+1}[-1]}\ \bigoplus\ \prod_{n\geq 1} \widehat{O}\RGraphs_{d+1}(1,n)\ot_{\bS_n} \sgn_n^{|d|}[d-dn]\\
    &\simeq & \OORGraphs_{d,d+1}(1,0)[d] \ (\text{isomorphism in the category of dg vector spaces}).
\Eeqrn   
The induced Lie bracket  in $\orgc_{d,d+1}$ is given via substitutions of the unique boundary of one generator of $\orgc_{d,d+1}$ into unlabelled {\it white}\, vertices of another generator.                                                                                           
There is a natural action of the dg Lie algebra $\orgc_{d+1}$ on $\orgc_{d,d+1}$
 given by substitutions of the 
unique boundary  of the generators of $\orgc_{d+1}$  into unlabelled 
{\it black}\, vertices of the generators of $\orgc_{d,d+1}$. Hence we can consider
 an auxiliary dg Lie algebra which, as a Lie algebra, is defined by the semidirect product,
\Beq\label{4: def of orgc_d,d+1}
\widehat{\orgc}_{d, d+1}:= \orgc_{d,d+1}\rtimes \orgc_{d+1},
\Eeq
and which is equipped with the twisted differential 
$$
d=\delta+\p, 
$$
where $\delta$ stands
for the original differential in the direct sum of complexes  $\orgc_{d,d+1}\oplus \orgc_{d+1}$ while
$\p$ is an injection
$$
\Ba{rccl}
\p: & \orgc_{d+1} & \lon & {\widehat{O}\RGraphs_{d+1}(1,0)[d]}\subset \orgc_{d,d+1}\\
& \Ga &\lon & (-1)^{|\Ga|} \Ga.
\Ea
$$

\subsubsection{\bf Remark} The dg Lie algebra $\widehat{\orgc}_{d, d+1}$
 can be equivalently defined 
as the deformation complex of a morphism of certain {\em 2-coloured}\,  properads,
$$
f: \Holie_{d,d+1} \lon \RGra_{d,d+1},
$$
which is {\it fully}\, analogous to the morphism 
$$
f: \Holie_{d,d+1} \lon \Gra_{d,d+1}
$$
studied in Proposition 3.3 of \cite{MWW}; the only difference is that one has to replace
the operad $\Gra_{d,d+1}$ in \cite{MWW} by its ribbon analogue $\RGra_{d,d+1}$; we omit the straightforward details.

\subsection{An interpolation between $\rgc_d$ and $\orgc_{d+1}$}
There is a morphism of dg Lie algebras 
$$
\pi_1: \widehat{\orgc}_{d, d+1} \lon \orgc_{d+1}
$$
given by the  projection of the semidirect product (\ref{4: def of orgc_d,d+1})
onto the second factor. The morphism $\pi_1$ is a quasi-isomorphism because its kernel --- the complex 
$\orgc_{d,d+1}$ ---   is isomorphic (up to a degree shift) to the complex $\OORGraphs_{d,d+1}(1,0)$ which is acyclic by 
Lemma {\ref{3: Lemma on acyclicity of twORGraphs(0,m)}}.

\sip

There is also a morphism of dg Lie algebras 
$$
\pi_2: \widehat{\orgc}_{d, d+1} \lon \rgc_{d}
$$
given by setting to zero every ribbon graph which has at least one non-special black vertex and then
making the identification (\ref{3: from special bu to dotted edges}). The kernel of this morphism is isomorphic as a graded vector (not as a complex)  to the following direct sum of graded vector spaces
$$
\ker\pi_2\simeq \OORGraphs_{d,d+1}^\bu (1,0)[d] \oplus \orgc_{d+1} \subset  \OORGraphs_{d,d+1} (1,0)[d] \oplus \orgc_{d+1} \simeq \widehat{\orgc}_{d, d+1} 
$$
Using the short exact sequence (\ref{3: shoert exact seq with Zm}) and the identification 
$\ORGC_{d+1}^{(1)}[-1]=\orgc_{d+1}$, one obtains an isomorphism in the category of graded vector spaces (see the previous section for the definition of the summands),
$$
\ker\pi_2\simeq Z_{d,d+1}^{(1)}[-1] \oplus \ \orgc_{d+1}[-1]\oplus  \orgc_{d+1}.
$$
One can pick up a spectral sequence of the complex $\Ker \pi_2$ (cf.\ proof of 
Proposition 4.2 in \cite{MWW}) whose initial page has the differential $\p+ d$,
where $\p$ acts only the last two summands above via the isomorphism
$$
\p:  \orgc_d\lon  \orgc_{d+1}[-1],
$$ 
and second differential $d$ acts only on the summand 
 $Z_{d,d+1}^{(1)}[-1]$ precisely
exactly as in the proof of Lemma {\ref{3: Lemma on H(ORGCd) and d,d+1}} where its acyclicity was established. We conclude that the complex $\ker\pi_2$ is acyclic. Hence the above morphism 
$\pi_2$ is quasi-isomorphism.

\sip

We obtain therefore a diagram of explicit quasi-isomorphisms of dg Lie algebras,
$$
\orgc_{d+1}    \stackrel{\pi_1}{\longleftarrow}  \widehat{\orgc}_{d,d+1}  \stackrel{\pi_2}{\lon}
\rgc_{d}
$$
which implies that the dg Lie algebras $\orgc_{d+1}$ and $\rgc_d$ are 
$\Lie_\infty$ quasi-isomorphic. Theorem {\ref{1: Main theorem 2}} is proven.

\sip

It is worth noting that one can assume without loss of generality that the dg Lie algebra $\orgc_d$ 
is generated by ribbon quivers with no so called  {\it passing}\, 
vertices, that is, the bivalent vertices with one incoming edge and 
one outgoing edge (cf.\ \cite{Wi1}).

\bip

\bip

{\Large
\section{\bf Computation of the cohomology of the derivation complex\\ of 
the dg properad $\PCY_d$}
}

\mip

\subsection{A complex of ribbon quivers with hairs} Our purpose in this subsection is 
to build a collection of $\bS_p^{op}\times \bS_q$-modules in the category of dg vector spaces,
$$
\PCYf_d=\{\PCYf_d(p,q)\}_{p\geq 1,q\geq 0},
$$
out of the twisted properad $\tw\ORGra_d$ or ribbon quivers discussed in \S {\ref{2: subsec on RGra^or and RGraphs^or}}.

\sip

We are interested only in the part
$\tw\ORGra_d(m,0)$
of the twisted dg properad $\tw\ORGra_d$ which is generated by ribbon quivers with no labelled vertices, and which therefore can not be composed with respect to the given properadic structure in  $\tw\ORGra_d$. Thus all vertices of generators $\Ga\in \tw\ORGra_d(m,0)$ are unlabelled and assigned the degree $d$.
The differential in $\tw\ORGra_d(m,0)$ is given by the general formula (\ref{2: delta in Tw(RGra^or)}) which in this case takes the form
\Beq\label{5: delta in Tw(RGra^or)(m,0)}
\delta\Ga:=
\sum_{i=1}^m \left(\hspace{-2mm} \Ba{c}\resizebox{4mm}{!}{  \xy
 (0,7)*{\bu}="A";
 (0,0)*+{_1}*\frm{o}="B";
 \ar @{->} "A";"B" <0pt>
\endxy} \Ea
+(-1)^d\hspace{-2mm}
\Ba{c}\resizebox{4mm}{!}{  \xy
 (0,7)*{\bu}="A";
 (0,0)*+{_1}*\frm{o}="B";
 \ar @{<-} "A";"B" <0pt>
\endxy} \Ea
\right)
\  _1\hspace{-0.7mm}\circ_i \Ga\ \
- \ \ (-1)^{|\Ga|}\hspace{-1mm}\sum_{v\in V_\bu(\Ga)} \Ga\circ_v  \left(\xy
 (0,0)*{\bu}="a",
(5,0)*{\bu}="b",
\ar @{->} "a";"b" <0pt>
\endxy\right).
\Eeq
This differential commutes with the action $\bS_m$ on $\tw\ORGra_d(m,0)$ given by relabelling of boundaries, and induces therefore a well-defined differential $\delta$ in the graded space
$$
\tw\ORGra_d(m,0)\ot_{\bS_m} \id_m
$$
which is generated by ribbon quivers with both vertices and boundaries unlabelled.
Let us call a vertex {\it source}\, (resp., {\it target})  if it has no incoming edges
(resp., no outgoing) edges.
The complex 
$$
C^{or}_d:= \prod_{m\geq 1}\left(\tw\ORGra_d(m,0)\ot_{\bS_m}\id_m, \delta\right)
$$
contains a subcomplex $I_d$ generated by ribbon quivers which have at least
one target of valency $\geq 2$. Hence we obtain
 a well defined quotient complex $Q_d:=C_d/I_d$ which is generated by ribbon quivers
$\Ga$ whose only targets are univalent. The set of vertices of each generator
$\Ga\in Q_d$ splits into the disjoint union
$$
V(\Ga)=V^1_{targets}(\Ga) \sqcup V^1_{sources}(\Ga) \sqcup V_{\geq 2}(\Ga)
$$
where $V^1_{targets}(\Ga)$ (resp. $V^1_{sources}(\Ga)$) is the set of univalent targets (resp., univalent sources), and  $V_{\geq 2}(\Ga)$ is the set of vertices of valency $\geq 2$. Note that each element  of $V_{\geq 2}(\Ga)$ has at least one outgoing edge.
As ribbon graphs are oriented, the set  $V^1_{targets}(\Ga)$ is non-empty.

\sip

The differential $\delta$ preserves the sets  $V^1_{targets}(\Ga)$ and $\sqcup V^1_{sources}(\Ga)$  $V^1_{in}(\Ga)$ of each generator  $\Ga\in  Q_d$ so that one obtains from $Q_d$  a well-defined
complex $Q_d^{marked}$ of  ribbon quivers $\Ga$  with {\it marked}\, univalent vertices, i.e.\ the ones which are equipped with some fixed isomorphisms of sets
$$
V^1_{targets}(\Ga) \rar [\# V^1_{targets}(\Ga)], \ \ \ \  V^1_{sources}(\Ga) \rar [\# V^1_{sources}(\Ga)].
$$
Let $Q_d^{marked}(p,q)$ be a subcomplex of $Q_d^{marked}$ generated by ribbon quivers
with $\# V^1_{sources}(\Ga)=q$ and $\# V^1_{targets}(\Ga)=p$; it is clearly an $\bS_p^{op}\times \bS_q$-module.
From now on we understand  a univalent target (resp.\, a univalent source) together with their 
unique attached directed edge  as an out-hair (resp.,  in-hair), 
and show in pictures correspondingly, e.g.

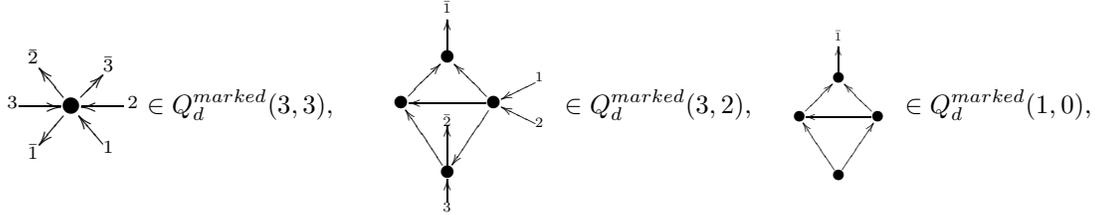
$$
\xy
(-5,-6)*{_{\bar{1}}}="a1",
(-5,6)*{^{\bar{2}}}="a2",
(5,5)*{^{\bar{3}}}="a3",
(5,-6)*{^{1}}="a4",
(8,0)*{\, ^2}="a5",
(-8,0)*{\, ^3}="a6",
(0,0)*{\bbu}="b",
\ar @{<-} "a1";"b" <0pt>
\ar @{<-} "a2";"b" <0pt>
\ar @{<-} "a3";"b" <0pt>
\ar @{->} "a4";"b" <0pt>
\ar @{->} "a5";"b" <0pt>
\ar @{->} "a6";"b" <0pt>
\endxy\in Q_d^{marked}(3,3), \ \ \ \
\Ba{c}\resizebox{21mm}{!}{
\mbox{$\xy
 (0,16)*{^{\bar{1}}}="out1";
 (0,-3)*{_{\bar{2}}}="out2";
 (0,-18)*{_3}="in3";
 (16,4)*{^{1}}="in1";
  (16,-4)*{^{2}}="in2";
 (-8,0)*{\bbu}="0";
 (0,8)*{\bbu}="1";
  (8,0)*{\bbu}="2";
(0,-12)*{\bbu}="3";
 \ar @{<-} "1";"0" <0pt>
  \ar @{<-} "1";"2" <0pt>
   \ar @{->} "2";"3" <0pt>
 \ar @{->} "2";"0" <0pt>
  \ar @{<-} "0";"3" <0pt>
  \ar @{->} "3";"out2" <0pt>
  \ar @{->} "1";"out1" <0pt> 
   \ar @{->} "in3";"3" <0pt>
    \ar @{<-} "2";"in1" <0pt> 
     \ar @{<-} "2";"in2" <0pt> 
\endxy$}}
\Ea
\in  Q_d^{marked}(3,2), \ \
\Ba{c}\resizebox{13mm}{!}{
\mbox{$\xy
 (0,16)*{^{\bar{1}}}="out1";
 (-8,0)*{\bbu}="0";
 (0,8)*{\bbu}="1";
  (8,0)*{\bbu}="2";
(0,-12)*{\bbu}="3";
 \ar @{<-} "1";"0" <0pt>
  \ar @{<-} "1";"2" <0pt>
   \ar @{<-} "2";"3" <0pt>
 \ar @{->} "2";"0" <0pt>
  \ar @{<-} "0";"3" <0pt>
  \ar @{->} "1";"out1" <0pt> 
\endxy$}}
\Ea
\in  Q_d^{marked}(1,0), 
$$
To fit this new terminology we change our notation
$$
H_{out}(\Ga):= V_{targets}^1(\Ga), \ \ \ H_{in}(\Ga):= V_{sources}^1(\Ga), \ \ 
V_{\geq 2}(\Ga)=:V(\Ga).
$$ 
and re-define the cohomological degree of such a ribbon quiver with hairs as follows,
$$
|\Ga|=d\#V(\Ga) + (1-d) \# E(\Ga) + (2-d)\# H_{out}(\Ga) - \# H_{in}(\Ga)
$$
 Every such a graph is tacitly assumed to be equipped with an orientation 
 $or(\Ga)$ which is, by definition, a unital vector in 
 \Bi
 \item[(i)] $(\det E(\Ga))\ot \det(H_{in}(\Ga))$  for $d$ even,
 \item[(ii)] $(\det V(\Ga))\ot \det(H_{out}(\Ga)\sqcup H_{in}(\Ga))$  for $d$ odd.
 \Ei

\subsection{A dg properad $\PCYf_d$} We make the $\bS$-bimodule $\{Q_d^{marked}(p,q)\}_{p\geq 1,q\geq 0}$ into a (non-unital) dg properad,
$$
\PCYf_d=\left\{\PCYf_d(p,q):= Q_d^{marked}(p,q) \right\}
$$
 whose composition is given naively by
gluing out-hairs of one ribbon quiver to in-hairs of another ribbon quiver creating 
thereby a new ribbon with hairs. The newly created internal edges get the degrees
$2-d-1=1-d$ so that such compositions have degree zero as required.
It follows from (\ref{5: delta in Tw(RGra^or)(m,0)}) that induced from 
$\tw\ORGra_d$ the
differential $\delta$  acts on a generator $\Ga$ of $\PCYf_d(p,q)$
by splitting the internal vertices
\Beq\label{5: delta in PreCYd}
\delta\Ga= (-1)^{|\Ga|}\hspace{-1mm}\sum_{v\in V_\bu(\Ga)} \Ga\circ_v  \left(\xy
 (0,0)*{\bu}="a",
(5,0)*{\bu}="b",
\ar @{->} "a";"b" <0pt>
\endxy\right)
\Eeq
in all possible ways which do not create a new univalent vertex. 
By valency of a vertex $v$ we understand the total number of edges and hairs attached.
 A {\it source} is a vertex with no incoming edges or in-hairs.
This differential is obviously consistent with the properadic structure so that we obtain
indeed a  dg properad $(\PCYf_d,\delta)$ of ribbon graphs with hairs.

\sip

The constructed dg properad  $(\PCYf_d,\delta)$ is obviously free as a properad.
It is generated by ribbon graphs which have precisely one vertex equipped with  
$p$ out-hairs and $q$-in hairs such that $p \geq 1$, $q\geq 0$. Let us call such 
ribbon graphs {\it haired ribbon $(p,q)$-corollas}, and let us denote 
the linear space generated by such ribbon corollas by $E_d(p,q)$. Each such corolla has the cohomological degree equal to $d+ (2-d)p - q$.
For example,
$$
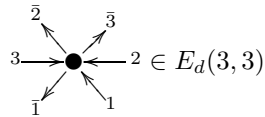
\xy
(-5,-6)*{_{\bar{1}}}="a1",
(-5,6)*{^{\bar{2}}}="a2",
(5,5)*{^{\bar{3}}}="a3",
(5,-6)*{^{1}}="a4",
(8,0)*{\, ^2}="a5",
(-8,0)*{\, ^3}="a6",
(0,0)*{\bbu}="b",
\ar @{<-} "a1";"b" <0pt>
\ar @{<-} "a2";"b" <0pt>
\ar @{<-} "a3";"b" <0pt>
\ar @{->} "a4";"b" <0pt>
\ar @{->} "a5";"b" <0pt>
\ar @{->} "a6";"b" <0pt>
\endxy\in E_d(3,3)
$$
If one sets $E_d:=\{E_d(p,q)\}_{p\geq 1, q\geq 0}$, then  
$$
\PCYf_d=\left(\cF ree \left\langle E_d\right\rangle, \delta \right)
$$
There is a one-to-one correspondence between representations of $\PCYf_d$ in a graded vector
space and degree $d$ pre-Calabi-Yau structures in $A$ \cite{LV}.

\subsection{A quotient properad $\PCY_d$} Let $I_s$ be the dg ideal in $\PCYf_d$
generated by haired ribbon graphs with at least one source or at least one passing vertex
(that is a bivalent vertex with precisely one incoming edge or hair, and precisely one 
outgoing edge or hair). The quotient properad
$$
\PCY_d:=\PCYf_d/I_s
$$
is spanned by ribbon quivers with all vertices at least trivalent, and with no sources
or targets. Every such a graph must have at least one in-hair and at least one out-hair.
The dg free properad $\PCY_d$ was first constructed in \cite{Q} as a cobar construction on the Koszul dual of the properad $\mathsf{BIB}_d$ of so called {\it balanced infinitesimal bialgebras}. For $d=2$ there is an epomorphism  \cite{Q}
$$
\PCY_2 \lon \mathsf{DPoiss}_\infty
$$
to the dg properad  $\mathsf{DPoiss}_\infty$ governing strongly homotopy double Poisson algebras 
which was described explicitly and studied  in \cite{LV}.

\sip

The dg properad $(\PCY_d,\delta)$ is uniquely characterized by the following property:
 {\em there is a one-to-one correspondence between 
representations of $\PCY_d$ in a dg vector space $(A,d)$ and  pre-Calabi-Yau 
structures in $A$}. Indeed, any a representation
$$
\rho: \PCY_d \lon \cE nd_A
$$
is uniquely determined by its values on the total vector space of generators,
$$
E_d^{total}:= \bigoplus_{p\geq 1,q\geq 1\atop p+q\geq 3} E_d(p,q),
$$
and it is easy to see that the graded vector space of equivariant linear maps
$$
\bigoplus_{p\geq 1,q\geq 1\atop p+q\geq 3} 
\Hom_{\bS_p^{op}\times \bS_p}\left(E_d(p.q), \cE nd_A(p,q)\right)
$$
can be identified  with the space 
  $$
{C}_{[d]}^3(A) := \prod_{k\geq 1}
\left( \bigoplus_{n_1,\ldots, n_k\geq 0\atop {n_1+...+n_k\geq 1\atop n_1+..+n_k+k\geq 3}}
\Hom\left(\bigotimes_{i=1}^{k}
(A[1])^{\ot n_i}, (A[2-d])^{\ot k}\right)_{\Z_k}\right).
$$
Thus every morphism $\rho: \PCY_d \lon \cE nd_A$
gives us a uniquely defined element $\pi\in {C}_{[d]}^3(A) $ of degree $d$, 
and the compatibility of $\rho$ with the differentials $\delta$ in $\PCY_d$ and $d$ in $A$
is equivalent to saying that $\pi$ satisfies the equation,
$$
d\pi +\frac{1}{2}[\pi,\pi]=0,
$$
where $[\ ,\ ]$ is the Lie bracket in ${C}_{[d]}^3(A)$ introduced and studied in \cite{IK,IKV,KTV}. 
Such MC elements give us a subclass of pre-CY algebra structures in $A$.

\subsection{The full derivation complex of $\PCY_d$} Let $\widehat{\PCY_d}$ be a completion of 
$\PCY_d$ with respect to the genus of the generating ribbon graphs. By a derivation complex of $\PCY_d$ we always understand (cf.\ \cite{MW2}) the derivation complex of its genus completion.
As $\PCY_d$ is a free properad, a derivation
$D$  acts on a haired ribbon graph $\Ga\in \PCY_d$ as a sum
$$
D(\Ga):=\sum_{v\in V(\Ga)} D(v)
$$
where the symbol $D(v)$ means the action of $D$ on the ribbon corolla $v$, that is, replacing $v$ with a suitable ribbon graph with hairs. Thus, as a graded vector space,
the space $\Der^{full}(\PCY_d)$  of all possible (genus completed) derivations is given by
$$
 \Der^{full}(\PCY_d)=
\prod_{p,q} \Hom_{\bS_p^{op}\times \bS_p}\left(E_d(p,q), \widehat{\PCY_d}(p,q)\right).
$$
It has the standard Lie bracket denoted by $[\ ,\ ]$.
The differential $\delta$ in $\PCY_d$ is an MC element of this Lie algebra, and hence it
makes  $\Der^{full}(\PCY_d)$ into a {\it dg}\, Lie algebra equipped with the differential given by
$$
d(D):=[\delta, D].
$$
The full complex $\Der^{full}(\PCY_d)$ is an awkward object in the combinatorial sense as hairs attached to corners of different boundaries 
have no {\em natural}\, cyclic ordering which the generators of $\PCY_d$ have. We study in this paper
its two subcomplexes whose actions on elements of $\PCY_d$ {\em preserve the number and the type of boundaries}\, of the generating ribbon graphs $\Ga$ from $\PCY_d$ and which admit a nice combinatorial description; thus the interpretation of generators of $\PCY_d$ as {\it ribbon}\, graphs with hairs plays a central role in what we dusciss next.

\subsection{The dg Lie algebra $\orgc_d$ as a subalgebra of $\Der^{full}(\PCY_d)$}
Let $\ga$ be any ribbon quiver from the dg Lie algebra $\orgc_d$ studied above. There is an associated derivation $D_\ga\in \Der^{full}(\PCY_d)$  which is defined by its action $D_\ga(v)$ on an arbitrary vertex $v$ of an
arbitrary element of $\PCY_d$ as follows: $D_\ga(v)$ means the substitution of the quiver $\ga$ into the vertex $v$ and taking the sum over attaching half-legs\footnote{A vertex $v$ of a haired ribbon graph can have edges and hairs of two (incoming and outgoing) types attached; we call all of them {\it half-legs attached to $v$}; the set of such half-legs has a natural cyclic ordering.} of $v$ to corners of the unique boundary of the quiver $\ga$  in all possible ways  while repsecting the cyclic orderings
of the set of half-legs and the set of corners of $\ga$ (this is a kind of properadic composition into $v$ as in the picture (\ref{2: composition in RGra}) above). It is easy
to see that this association $\ga\rar D_\ga$ gives us a monomorphism
of dg Lie algebras
$$
0\lon \orgc_d \lon \Der^{full}(\PCY_d)
$$
whose image is denoted by the same symbol $\orgc_d$. It is worth noting that the differential $\delta$ in $\PCY_d$ acts as the derivation $D_{\tau}$ where
\Beq\label{5: tau}
\tau:=
{\xy
 (0,0)*{\bu}="a",
(5,0)*{\bu}="b",
\ar @{->} "a";"b" <0pt>
\endxy}.
\Eeq

\subsection{A dg Lie algebra $\Der(\PCY_d)$ of ``preserving boundaries" derivations}
\label{5: subsec on Der(PCY)}
Let $\PCY_d(1;p,q)$ be  the linear subspace of $\PCY_d(p,q)$ generated by haired ribbon graphs with
precisely one boundary, and let 
$\widehat{\PCY_d}(1;p\oplus q)$ be a version of
 $\widehat{\PCY_d}(1;p,q)$ in which the numerical labels of all hairs are forgotten, i.e.\ the set of all hairs attached to a ribbon quiver $\Theta$ from $\PCY_d(1;p\oplus q)$ is equipped only with the induced (from the set of corners of the unique boundary) cyclic order. If we shrink to a point all edges of $\Theta$ we obtain an unlabelled ribbon corolla $v$ of with  $p+q$ hairs; we say that such a ribbon corolla $v$
 is {\em compatible}\, with $\Theta$. We claim that the graded vector space
 \Beq\label{5: Der as product}
\Der(\PCY_d)\simeq \prod_{p,q} \widehat{\PCY_d}(1;p\oplus q)[q-(2-d)p-d]
\Eeq
is a dg Lie subalgebra of  $\Der^{full}(\PCY_d)$. Indeed, for any generator 
$\Theta\in  \Der(\PCY_d)$ there is an associated derivation $D_\Theta$ of the properad $\PCY_d$
which acts on a vertex $v$ of any generator $\Ga\in \PCY_d$ as follows,
\Bi
\item[(i)] we set $D_\Theta(v)=0$ if $v$ is not compatible with $\Theta$, otherwise
\item[(ii)] we set $D_\Theta(v)$ to be the sum  over all possible gluings of half-legs  of the ribbon corolla $v$ to hairs of the quiver $\Theta$  while respecting cyclic orderings of both sets and their ("in" or "out") types. 
\Ei
It is easy to see that the subspace
$$
\Der(\PCY_d)\subset \Der^{full}(\PCY_d)
$$
is a Lie subalgebra. Moreover, it is a {\it dg}\, Lie subalgebra as the differential $\delta$ in $\PCY_d$ can be identified with the derivation $D_\Delta$ associated to the following series of unlabelled haired ribbon graphs,
$$
\Delta:=\sum_{p,q\geq 1}  \underbrace{
\sum_{\text{maps}\ [p+q]\rar C(b)} \resizebox{14mm}{!}{\xy
(-5,-6)*{}="a1",
(-5,6)*{}="a2",
(5,5)*{}="a3",
(5,-6)*{}="a4",
(8,0)*{}="a5",
(-8,0)*{}="a6",
(0,0)*+{\tau}="b",
\ar @{<-} "a1";"b" <0pt>
\ar @{<-} "a2";"b" <0pt>
\ar @{<-} "a3";"b" <0pt>
\ar @{->} "a4";"b" <0pt>
\ar @{->} "a5";"b" <0pt>
\ar @{->} "a6";"b" <0pt>
\endxy}}_{\text{ribbon $(p,q)$-corolla}}\ \in \Der(\PCY_d)
$$
where $\tau$ is given in (\ref{5: tau}) and the sum is taken over all possible cyclically inequivalent ways to attach  $p+q$ hairs to the two corners of the
 unique boundary $b$ of $\tau$, and setting to zero every resulting summand which has at least one vertex of 
valency $\leq 2$, or at least one target, or at least one source. Hence
the differential $d$ in  $\Der(\PCY_d)$ consists of two terms
\Beq\label{5: d in Der(PCY)}
d=\delta + \delta'
\Eeq
where $\delta$ is the standard differential (\ref{5: delta in PreCYd}) in 
$\PCY_d$ and $\delta'$ acts on a generator
$\Ga\in \PCY_d(1;p\oplus q)$ by attaching haired ribbon corollas
to every in-hair of $\Ga$ and any out-hair of $\Ga$ as shown schematically in 
the following picture,
$$
\delta'\Ga=\sum_{p,q\geq 1\atop p+q\geq 3}\hspace{-5mm}
\underbrace{\resizebox{12mm}{!}{
 \xy
(-5,-16)*{}="a1",
(-5,-4)*{}="a2",
(6,0)*{\Ga}="a3",
(5,-16)*{}="a4",
(8,-10)*{}="a5",
(-8,-10)*{}="a6",
(0,-10)*{\bbu}="b",
\ar @{<-} "a1";"b" <0pt>
\ar @{<-} "a2";"b" <0pt>
\ar @{<-} "a3";"b" <0pt>
\ar @{->} "a4";"b" <0pt>
\ar @{->} "a5";"b" <0pt>
\ar @{->} "a6";"b" <0pt>
\endxy}}_{\text{ribbon $(p,q)$-corolla}}
\pm
\sum_{p,q\geq 1\atop p+q}\hspace{-5mm}
\overbrace{\resizebox{12mm}{!}{
 \xy
(-5,0)*{}="a1",
(-5,12)*{}="a2",
(6,12)*{}="a3",
(5,-4)*{\Ga}="a4",
(8,6)*{}="a5",
(-8,6)*{}="a6",
(0,6)*{\bbu}="b",
\ar @{<-} "a1";"b" <0pt>
\ar @{<-} "a2";"b" <0pt>
\ar @{<-} "a3";"b" <0pt>
\ar @{->} "a4";"b" <0pt>
\ar @{->} "a5";"b" <0pt>
\ar @{->} "a6";"b" <0pt>
\endxy}}^{\text{ribbon $(p,q)$-corolla}}
$$

\subsubsection{\bf Rescaling cohomology class in $H^\bu(\Der(\PCY_d))$}
The properad $\PCY_d$ admits an obvious automorphism $r_\lambda$ given on the generating
ribbon $(p,q)$-corollas $c_{p.q}$ by a rescaling
$$
r_{\lambda}: c_{p,q} \lon \la^{p-q} c_{p,q}\ \ \forall\ \la\in \K\setminus 0.
$$
Differentiaing over $\lambda$ and setting $\la=1$ one obtains a non-trivial derivation of
$\PCY_d$ which is given by the following cohomology class in $H^\bu(\Der(\PCY_d))$
$$
r:=\sum_{p,q\geq 1\atop p+1\geq 3}(p-q)\hspace{-5mm}  \underbrace{
\resizebox{10mm}{!}{\xy
(-5,-6)*{}="a1",
(-5,6)*{}="a2",
(5,5)*{}="a3",
(5,-6)*{}="a4",
(8,0)*{}="a5",
(-8,0)*{}="a6",
(0,0)*{\bbu}="b",
\ar @{<-} "a1";"b" <0pt>
\ar @{<-} "a2";"b" <0pt>
\ar @{<-} "a3";"b" <0pt>
\ar @{->} "a4";"b" <0pt>
\ar @{->} "a5";"b" <0pt>
\ar @{->} "a6";"b" <0pt>
\endxy}}_{\text{ribbon $(p,q)$-corolla}}
$$
called the {\em rescaling class}; here the sum is taken  over all possible 
cyclically inequivalent ways to attached 
$p$ in-hairs and $q$-out hairs to one vertex.  This infinite sum  is obviously closed under $d$ as the ribbon
graphs originating from the splitting of the black vertex via $\delta r$ 
cancel out the  attachments terms coming from $\delta' r$. It can not be a cobounadry 
as $r$ is spanned by haired ribbon graphs with only one vertex.

\subsection{From $\orgc_d$ to $\Der(\PCY_d)$}
 There is a  morphism
of dg Lie algebras,
\Beq\label{3: Morhism F from dcGC to Der^++}
\Ba{rccc}
 F: & \orgc_d &\to & \Der(\PCY_d)\\
         &   \ga & \to & F(\Ga)
         \Ea
\Eeq
where the derivation $F(\ga)$ acts (from the right) on the generating ribbon
$(p,q)$-corolla of $\PCY_b$ by attaching  $p+q$ hairs to the corners of the
 unique boundary $b$ of $\Ga$ in all possible ways while respecting the cyclic 
 orders of both sets,
 $$
 F(\Ga): \underbrace{
\resizebox{10mm}{!}{\xy
(-5,-6)*{}="a1",
(-5,6)*{}="a2",
(5,5)*{}="a3",
(5,-6)*{}="a4",
(8,0)*{}="a5",
(-8,0)*{}="a6",
(0,0)*{\bbu}="b",
\ar @{<-} "a1";"b" <0pt>
\ar @{<-} "a2";"b" <0pt>
\ar @{<-} "a3";"b" <0pt>
\ar @{->} "a4";"b" <0pt>
\ar @{->} "a5";"b" <0pt>
\ar @{->} "a6";"b" <0pt>
\endxy}}_{\text{ribbon $(p,q)$-corolla}}
\ \ \ \lon \ 
 \sum_{\text{maps}\ [p+q]\rar C(b)} \resizebox{14mm}{!}{\xy
(-5,-6)*{}="a1",
(-5,6)*{}="a2",
(5,5)*{}="a3",
(5,-6)*{}="a4",
(8,0)*{}="a5",
(-8,0)*{}="a6",
(0,0)*+{\ga}="b",
\ar @{<-} "a1";"b" <0pt>
\ar @{<-} "a2";"b" <0pt>
\ar @{<-} "a3";"b" <0pt>
\ar @{->} "a4";"b" <0pt>
\ar @{->} "a5";"b" <0pt>
\ar @{->} "a6";"b" <0pt>
\endxy},
 $$
and setting to zero every graph in the r.h.s.\  which has at least one vertex of 
valency $\leq 2$, or at least one target, or at least one source. It is easy to see
that the map $F$ respects the differentials and Lie brackets (cf.\ \cite{MW2}).

\sip

Under isomorphism (\ref{5: Der as product}), the derivation $F(\ga)$ can be obtained
from any given ribbon quiver $\ga\in \orgc_d$ by taking an infinite sum over all possible  ways to attach
to {\em each} corner $c_i$ of {\it each}\, vertex of $v\in V(\ga)$ (whose attached  edges are shown now schematically as undirected dotted ones)  $p_i\geq 0$ in-hairs and $q_i$ out-hairs
for all possible values of $p_i$ and $q_i$ and in all possible orders.
$$
F:  v= \Ba{c}\resizebox{21mm}{!}{\xy
(4,2.1)*{{c_i}};
(0,3.8)*{{c_{j}}};
%
 (0,0)*{\xycircle(6,6){.}};
(-6,0)*{}="1";
(-16,0)*{}="1'";
(4.2,4.2)*{}="2";
(12,12)*{}="2'";
(4.2,-4.2)*{}="3";
(12,-12)*{}="3'";
(6,0)*{}="4";
(16,0)*{}="4'";
(-4.2,4.2)*{}="5";
(-12,12)*{}="5'";
(-4.2,-4.2)*{}="6";
(-12,-12)*{}="6'";
(-0.3,6)*{}="a";
(-0.4,6)*{}="b";
\ar @{.} "1";"1'" <0pt>
\ar @{.} "2";"2'" <0pt>
\ar @{.} "3";"3'" <0pt>
\ar @{.} "4";"4'" <0pt>
\ar @{.} "5";"5'" <0pt>
\ar @{.} "a";"b" <0pt>
\ar @{.} "6";"6'" <0pt>
\endxy}\Ea
\lon
\sum
\Ba{c}\resizebox{29mm}{!}{\xy
(0,12)*{{p_j+q_j}};
(0,8)*{{...}};
(-8,4)*{{...}};
(24,8)*{{p_i+q_i}};
(4,2.1)*{{c_i}};
(0,3.8)*{{c_{j}}};
(5,3.5)*{}="A1";
(14,9.5)*{}="B1";
(5.5,2.5)*{}="A2";
(16,6)*{}="B2";
(5.5,1.5)*{}="A3";
(18,3)*{}="B3";
 (0,0)*{\xycircle(6,6){.}};
(-6,0)*{}="1";
(-16,0)*{}="1'";
(4.2,4.2)*{}="2";
(12,12)*{}="2'";
(4.2,-4.2)*{}="3";
(12,-12)*{}="3'";
(6,0)*{}="4";
(16,0)*{}="4'";
(-4.2,4.2)*{}="5";
(-12,12)*{}="5'";
(-4.2,-4.2)*{}="6";
(-12,-12)*{}="6'";
(-0.3,6)*{}="a";
(-0.4,6)*{}="b";
\ar @{.} "1";"1'" <0pt>
\ar @{.} "2";"2'" <0pt>
\ar @{.} "3";"3'" <0pt>
\ar @{.} "4";"4'" <0pt>
\ar @{.} "5";"5'" <0pt>
\ar @{.} "a";"b" <0pt>
\ar @{.} "6";"6'" <0pt>
\ar @{->} "A1";"B1" <0pt>
\ar @{<-} "A2";"B2" <0pt>
\ar @{->} "A3";"B3" <0pt>
\endxy}\Ea
$$

If $v$ is a target in $\ga$, then  at least one out-hair must be attached to some its corner,
if $v$ is a source, then at least one in-hair must be attached to some its corner.

\sip

The ``vertex"  $F(v)$ is called the {\it full hairy}\, version of $v$. Therefore, the map $F$
becomes, in a sense, tautological: $F(\ga)$ is obtained from $\ga$ be declaring every vertex of $\ga$ full hairy. It is easy to see that this map $F$ has degree zero, and respects both the differentials and Lie brackets (cf.\ \cite{MW2}).

\subsubsection{\bf Remark}\label{5: remaik on full hairs}
For future reference, consider a version $F^{c_i}_{in}(v)$ (resp. $F^{c_i}_{out}(v)$) of the above construction of $F(v)$ in which we add only in-hairs  (resp. only out-hairs) to a particular corner $c_i$ of the vertex $v$, 
i.e. the summation  as above goes only over $q_i\geq 0$, $p_i=0$ (resp. over $p_i\geq 0$, $q_i=0$), and let us call
the formal sum $F^{c_i}_{in}(v)$ (resp., $F^{c_i}_{out}(v)$) the vertex $v$ with {\it  full in-hairy}\,  (resp. {\it full out-hairy}) corner $c_i$. If each corner of $v$ is full in-hairy
(resp., full out-hairy), then we call $v$ {\it full in-hairy}\,  (resp. {\it full out-hairy}).

\subsection{Theorem}\label{5: Theorem on F qis} {\it The morphism $F$ is a quasi-isomorphism up to 
one rescaling class $r$.} 
 
 \begin{proof} The proof follows the same ``trimming of branches" scenario which was used in 
 \cite{MW2} to prove a similar statement for the dg Lie algebra of oriented graphs $\OGC_{c+d+1}$ and the derivation complex of the dg properad $\HoLB_{c,d}$ of degree shifted
 strongly homotopy Lie bialgebras.

 \sip

 Call a vertex $v$ of a ribbon quiver $\Ga\in \Der(\PCY_d)$ {\it in-special}\ if it has no outgoing hairs,  precisely one outgoing edge, and any number of ingoing hairs or ingoing edges, i.e. $v$ looks like a generator of the dg operad $\cA ss_\infty$ controlling strongly homotopy associative algebras,
 $$
 v=  
 \Ba{c}\resizebox{22mm}{!}{ \xy
(1,-5)*{\ldots},
(-13,-8)*{},
(-8,-8)*{},
(-3,-8)*{},
(7,-8)*{},
(13,-8)*{},
 (0,0)*{\bu}="a",
(0,7)*{}="0",
(-12,-5)*{}="b_1",
(-8,-5)*{}="b_2",
(-3,-5)*{}="b_3",
(8,-5)*{}="b_4",
(12,-5)*{}="b_5",
\ar @{->} "a";"0" <0pt>
\ar @{<-} "a";"b_2" <0pt>
\ar @{<-} "a";"b_3" <0pt>
\ar @{<-} "a";"b_4" <0pt>
\endxy}\Ea.
 $$
 Call such a vertex {\it quasi-univalent}\,
 it has no incoming edges, only incoming hairs. Let $\Ga^{sk}$ be the graph 
 obtained from a generator $\Ga\in \Der(\PCY_d)$ by erasing iteratively all quasi-univalent in-special vertices; we call $\Ga^{sk}$ the {\it skeleton of} $\Ga$, and call its vertices {\it root or skeleton vertices}. Thus $\Ga$ is obtained
 from $\Ga^{sk}$ by attaching to the in-hairs (if any) of the root vertices $r\in V(\Ga^{sk})$ {\it in-trees}\, $T$ which can be identified with the elements of the  operad of 
  $\cA ss _\infty$ controlling strongly homotopy associative algebras,
 $$
T=\xy
(0,16)*+{R}*\frm{o}="R"; 
(-7,14)*{}="1";
(7,14)*{}="1'";
(-7,20)*{}="1''";
(4.2,4.2)*{}="2";
(12,12)*{}="2'";
(4,11)*{}="r1";
(-4,11)*{}="r2";
(-10.5,-2)*{},
(-11,-17)*{},
(-2,-17)*{},
(3,-10)*{},
(8,-10)*{},
(14,-10)*{},
(21,-10)*{},
(0,14)*{}="0",
 (0,8)*{\bu}="a",
(-10,0)*{}="b_1",
(-2,0)*{\bu}="b_2",
(12,0)*{\bu}="b_3",
(2,-8)*{}="c_1",
(-7,-8)*{\bu}="c_2",
(8,-8)*{}="c_3",
(14,-8)*{}="c_4",
(20,-8)*{}="c_5",
(-11,-15)*{}="d_1",
(-3,-15)*{}="d_2",
\ar @{->} "a";"R" <0pt>
\ar @{<-} "a";"b_1" <0pt>
\ar @{<-} "a";"b_2" <0pt>
\ar @{<-} "a";"b_3" <0pt>
\ar @{<-} "b_2";"c_1" <0pt>
\ar @{<-} "b_2";"c_2" <0pt>
\ar @{<-} "b_3";"c_3" <0pt>
\ar @{<-} "b_3";"c_4" <0pt>
\ar @{<-} "b_3";"c_5" <0pt>
\ar @{<-} "c_2";"d_1" <0pt>
\ar @{<-} "c_2";"d_2" <0pt>
\ar @{.} "R";"1" <0pt>
\ar @{.} "R";"1'" <0pt>
\ar @{.} "R";"1''" <0pt>
\ar @{<-} "R";"r1" <0pt>
\ar @{<-} "R";"r2" <0pt>
\endxy
$$
 The edges of the root vertex $R$ which belong to $\Ga^{sk}$ or its {\it out}-hairs  are shown schematically as dotted undirected half-edges; these dotted half-edges partition $r$ into so called {\it skeleton corners}. We define the valency $|c|$ of 
 such a skeleton corner $c$ as follows,
 $$
 |c|= \text{number of in-hairs attached to $c$} + \text{number of in-trees attached to $c$}.
 $$ 
 A corner $c$ of $R$ is called {\it bold}, if $|c|=0$, and {\it univalent}\, if $|c|=1$;
 otherwise $c$ is called {\it generic}.

 \sip

Consider a filtration of the complex $\Der(\PCY_d)$ by the number of vertices which are {\it not}\, in-special,
and let
 $\{E_r, d_r\}_{r\geq 0}$ be the associated spectral sequence. The induced differential $d_0=\delta_0 + \delta_0'$ acts on a generator $\Ga$ of the initial page $E_0$ 
 by attaching an  in-special vertex
 to in-hairs via $\delta_0'$  or by splitting vertices in such a way that at least of the newly created
 vertices is in-special.

\sip

 Let $C_{\text{full in-hairy}}$ be a subspace of $E_0$ generated by linear 
 combinations of ribbon quivers which have no in-special vertices and whose all 
 vertices are the {\it full in-hairy}\, as defined
in Remark {\ref{5: remaik on full hairs}}. This subspace is a trivial subcomplex
of $(E_0,d_0)$ as $\delta_0\Ga=- \delta_0'\Ga$ for any 
$\Ga\in C_{\text{full in-hairy}}$.

\sip

{\bf Claim A}. {\em The inclusion $C_{\text{full in-hairy}} \hook E_0$ is a quasi-isomorphism.}

\sip
We prove this claim in two steps.

\sip

{\sc Step 1}. Let us consider a filtration of $(E_0,d_0)$ by the number 
of skeleton vertices in the generators, and let $\{\cE_r(E_0),d_{0r}\}_{r\geq 0}$ stand for the associated spectral sequence. The induced
differential $d_{00}$ in $\cE_0(E_0)$ is the sum 
 $d_{00}=\delta_{00} + \delta_0'$,
where $\delta_0'$ is the same ``attaching in-special vertices" differential as in $E_0$, while $\delta_{00}$ is the usual ``splitting of vertices" differential  which acts on a corner $c$ with $|c|\geq 2$ of each root vertex $R\in V(\Ga^{sk})$ as well as on each vertex of any in-tree $T$ attached to $c$  exactly as in the dg operad $\cA ss_\infty$,
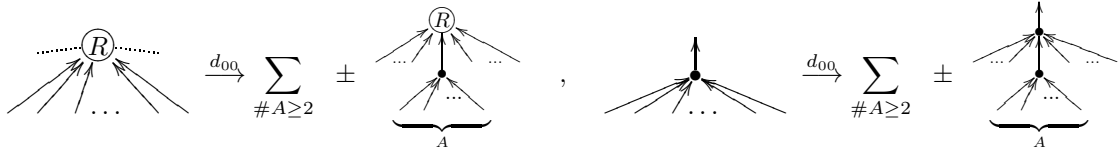
$$
\xy
(1,-5)*{\ldots},
(-13,-7)*{},
(-8,-7)*{},
(-3,-7)*{},
(7,-7)*{},
(13,-7)*{},
(0,4)*+{R}*\frm{o}="a"; 
(-8,3)*{}="1";
(8,3)*{}="1'";
(-12,-5)*{}="b_1",
(-8,-5)*{}="b_2",
(-3,-5)*{}="b_3",
(8,-5)*{}="b_4",
(12,-5)*{}="b_5",
%
\ar @{<-} "a";"b_2" <0pt>
\ar @{<-} "a";"b_3" <0pt>
\ar @{<-} "a";"b_1" <0pt>
\ar @{<-} "a";"b_4" <0pt>
\ar @{<-} "a";"b_5" <0pt>
\ar @{.} "a";"1" <0pt>
\ar @{.} "a";"1'" <0pt>
\endxy
\stackrel{d_{00}}{\lon}
\sum_{\# A\geq 2} \ \ \pm
 \Ba{c}\resizebox{20mm}{!}{  \xy
 (0,11)*{}="1";
  (-5,1)*{}="l'";
(0,8)*+{R}*\frm{o}="0"; 
(-8,7)*{}="1";
(8,7)*{}="1'";  
(5,1)*{}="r1";
(13,1)*{}="r2";
(-11,1)*{}="l2";
(9,1)*{_{...}};
(-7,1)*{_{...}};
%
(0,-1)*{\bu}="C";
(-7,-7)*{}="L1";
(-3,-7)*{}="L2";
(2,-5)*{...};
(7,-7)*{}="L3";
<0mm,-11mm>*{\underbrace{ \ \ \ \ \ \ \ \ \ \ \ \ \ \ }_{A}},
\ar @{<-} "C";"L1" <0pt>
\ar @{<-} "C";"L2" <0pt>
\ar @{<-} "C";"L3" <0pt>
\ar @{->} "C";"0" <0pt>
\ar @{->} "l'";"0" <0pt>
\ar @{->} "r1";"0" <0pt>
\ar @{->} "r2";"0" <0pt>
\ar @{->} "l2";"0" <0pt>
 \endxy}
 \Ea
 \ \ \ , \ \ \ 
\xy
(1,-5)*{\ldots},
(-13,-7)*{},
(-8,-7)*{},
(-3,-7)*{},
(7,-7)*{},
(13,-7)*{},
 (0,0)*{\bu}="a",
(0,5)*{}="0",
(-12,-5)*{}="b_1",
(-8,-5)*{}="b_2",
(-3,-5)*{}="b_3",
(8,-5)*{}="b_4",
(12,-5)*{}="b_5",
\ar @{->} "a";"0" <0pt>
\ar @{<-} "a";"b_2" <0pt>
\ar @{<-} "a";"b_3" <0pt>
\ar @{<-} "a";"b_1" <0pt>
\ar @{<-} "a";"b_4" <0pt>
\ar @{<-} "a";"b_5" <0pt>
\endxy
\stackrel{d_{00}}{\lon}
\sum_{\# A\geq 2} \ \ \pm
 \Ba{c}\resizebox{20mm}{!}{  \xy
 (0,11)*{}="1";
  (-5,1)*{}="l'";
(-0,6)*{\bu}="0";
(5,1)*{}="r1";
(13,1)*{}="r2";
(-11,1)*{}="l2";
(9,1)*{_{...}};
(-7,1)*{_{...}};
%
(0,-1)*{\bu}="C";
(-7,-7)*{}="L1";
(-3,-7)*{}="L2";
(2,-5)*{...};
(7,-7)*{}="L3";
<0mm,-11mm>*{\underbrace{ \ \ \ \ \ \ \ \ \ \ \ \ \ \ }_{A}},
\ar @{<-} "C";"L1" <0pt>
\ar @{<-} "C";"L2" <0pt>
\ar @{<-} "C";"L3" <0pt>
\ar @{->} "C";"0" <0pt>
\ar @{<-} "1";"0" <0pt>
\ar @{->} "l'";"0" <0pt>
\ar @{->} "r1";"0" <0pt>
\ar @{->} "r2";"0" <0pt>
\ar @{->} "l2";"0" <0pt>
 \endxy}
 \Ea
$$
where the summation goes over {\it connected}\, subsets $A$ of the original totally ordered set of
all ingoing edges/hairs. 

\sip 

Let $B_{\text{corners}}$ be a subspace (in fact, a trivial subcomplex)
of $\cE_0(E_0)$ generated by linear combinations of ribbon quivers which have no in-special vertices, and whose vertices have either bold corners or full hairy corners. 
Such ribbon quivers are cycles with respect to the differential $d_{00}$. 

\sip

{\bf Claim A(i)}. {\em The inclusion $B_{\text{corners}} \hook \cE_{0}(E_0)$ is a quasi-isomorphism.}

\sip

Indeed, let us consider a filtration of  $\cE_{0}(E_0)$ by the total number of in-hairs
of its generators. The induced differential in the associated graded $gr \cE_{0}(E_0)$
is precisely the $A_\infty$-like differential $d_{00}$. We conclude almost immediately
that $H^\bu(gr \cE_{0}(E_0))$ is generated by skeleton-like graphs $\Ga$ whose vertices have either bold corners  or corners with precisely one incoming hair,
$$
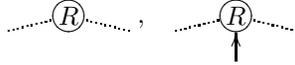
\xy
(0,0)*+{R}*\frm{o}="R"; 
(-8,-2)*{}="1";
(8,-2)*{}="1'";  
\ar @{.} "R";"1" <0pt>
\ar @{.} "R";"1'" <0pt>
\endxy
\ , \ \ \ 
\xy
(0,0)*+{R}*\frm{o}="R"; 
(-8,-2)*{}="1";
(8,-2)*{}="1'";  
(0,-6)*{}="r1";
\ar @{.} "R";"1" <0pt>
\ar @{.} "R";"1'" <0pt>
\ar @{->} "r1";"R" <0pt>
\endxy
$$

 This space is identical to $B_{\text{corners}}$. The Comparison of Spectral Sequences  Theorem completes the proof of the {\bf Claim A(i)}.

\sip

{\sc Step 2}. We conclude from Step 1 that the next page $(\cE_1(E_0), d_{10})$
can be identified  $B_{\text{corners}}$, and the induced differential $d_{10}$ acts
on the generators $\Ga\in B_{\text{corners}}$ by splitting skeleton vertices in such a way
that at least one of the newly created vertices is a in-special vertex of the form
$$
\xy
(0,0)*+{R}*\frm{o}="R"; 
(-9,0)*{}="1";
(9,0)*{}="2";  
\ar @{<.} "R";"1" <0pt>
\ar @{.>} "R";"2" <0pt>
\endxy
$$
i.e. it has one incoming skeleton edge, one outgoing skeleton edge and precisely two corners such that at least one of these corners is full in-hairy. Let us call such skeleton 
vertices {\it passing}.

\sip

The subspace $C_{\text{full in-hairy}}$ considered in {\bf Claim A} is a trivial subcomplex of $B_{\text{corners}}$ generated by skeleton graphs with no passing vertices and whose every vertex is full in-hairy (as defined in Remark  {\ref{5: remaik on full hairs}}).

\sip

{\bf Claim A(ii)}. {\em The inclusion 
$$
C_{\text{\rm full in-hairy}} \hook \cE_{1}(E_0)\simeq B_{\text{\rm corners}}
$$ 
is a quasi-isomorphism.}

\sip

Indeed, as every generator $\Ga$ of $B_{\text{corners}}$ is oriented, it must have at least
one vertex with no incoming skeleton edges. Hence $\Ga$ has at least one corner which is 
full in-hairy. The {\bf Claim A(ii)} is proven if we show that the subcomplex $B'$ of $B_{\text{corners}}$ generated by skeleton graphs $\Ga$ with at least one vertex having a bold corner $c$ is acyclic. Since $\Ga$ has only one boundary $b$, there is a pair of 
corners $c_i$ and $c_j$ in $b$ such that
\Bi
\item[(i)] the corner $c_i$ is bold while the corner $c_j$ is full hairy,
\item[(ii)] the corners $c_i$ and $c_j$ do not belong to passing vertices,
 \item[(iii)] the corners $c_i$ and $c_j$ are connected to each other in the boundary $b$
 either by one edge or by a directed path which goes only through the fully hairs corners $c_{passing}$   of some passing vertices of $\Ga$.
\Ei
A standard argument based on the number of fully hairy passing corners $c_{passing}$ on the directed path between
$c_i$ to $c_j$ (cf.\ \cite{Wi1,MW2})
shows that $B'$ is acyclic. The  {\bf Claim 2(ii)} is proven. Hence the {\bf Claim A} follows.

\sip

{\sc Step 3}. We study next the page $E_1\simeq C_{\text{full in-hairy}}$ of the spectral sequence by the number of in-special vertices. The induced differential $d_1$
creates now vertices with $\geq 2$ out-going edges or hairs. One  repeats essentially  Steps 1 and 2
 by considering a filtration by  the number of so called
{\it out-special vertices}\, $v$, 
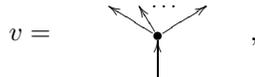
$$
 v=  
 \Ba{c}\resizebox{22mm}{!}{ \xy
(1,5)*{\ldots},
(-13,8)*{},
(-8,8)*{},
(-3,8)*{},
(7,8)*{},
(13,8)*{},
 (0,0)*{\bu}="a",
(0,-7)*{}="0",
(-12,5)*{}="b_1",
(-8,5)*{}="b_2",
(-3,5)*{}="b_3",
(8,5)*{}="b_4",
(12,5)*{}="b_5",
\ar @{<-} "a";"0" <0pt>
\ar @{->} "a";"b_2" <0pt>
\ar @{->} "a";"b_3" <0pt>
\ar @{->} "a";"b_4" <0pt>
\endxy}\Ea,
 $$
which are defined simply by reversing directions in the definition of in-special vertices. 
One introduces similarly the notion of the skeleton $\Ga^{sk}$ of any generator $\Ga\in E_1$, then studies the
complex geneerated by rooted trees of the form
  $$
T=
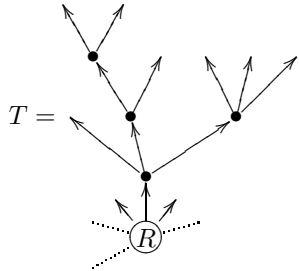
\xy
(0,-16)*+{R}*\frm{o}="R"; 
(-7,-14)*{}="1";
(7,-14)*{}="1'";
(-7,-20)*{}="1''";
(4.2,-4.2)*{}="2";
(12,-12)*{}="2'";
(4,-11)*{}="r1";
(-4,-11)*{}="r2";
(-10.5,2)*{},
(-11,17)*{},
(-2,17)*{},
(3,10)*{},
(8,10)*{},
(14,10)*{},
(21,10)*{},
(0,14)*{}="0",
 (0,-8)*{\bu}="a",
(-10,0)*{}="b_1",
(-2,0)*{\bu}="b_2",
(12,0)*{\bu}="b_3",
(2,8)*{}="c_1",
(-7,8)*{\bu}="c_2",
(8,8)*{}="c_3",
(14,8)*{}="c_4",
(20,8)*{}="c_5",
(-11,15)*{}="d_1",
(-3,15)*{}="d_2",
\ar @{<-} "a";"R" <0pt>
\ar @{->} "a";"b_1" <0pt>
\ar @{->} "a";"b_2" <0pt>
\ar @{->} "a";"b_3" <0pt>
\ar @{->} "b_2";"c_1" <0pt>
\ar @{->} "b_2";"c_2" <0pt>
\ar @{->} "b_3";"c_3" <0pt>
\ar @{->} "b_3";"c_4" <0pt>
\ar @{->} "b_3";"c_5" <0pt>
\ar @{->} "c_2";"d_1" <0pt>
\ar @{->} "c_2";"d_2" <0pt>
\ar @{.} "R";"1" <0pt>
\ar @{.} "R";"1'" <0pt>
\ar @{.} "R";"1''" <0pt>
\ar @{->} "R";"r1" <0pt>
\ar @{->} "R";"r2" <0pt>
\endxy
$$
and arrives --- in the full analogy to Steps 1 and 2 above --- to the conclusion that 
the inclusion $F(\orgc_d)\hook E_1$ is a quasi-isomorphism. We do not repeat the arguments shown 
in Steps 1 and 2 just above.
 \end{proof}
 
 Theorem {\ref{1: Theorem on H(Der)}} in the Introduction follows immediately from Theorem 
 {\ref{5: Theorem on F qis}},
  Theorem {\ref{1: Main theorem 1}} and Theorem {\ref{1: Main theorem 2}}.
 
 \subsection{Some applications} It is well-known that $H^k_c(\cM_{g,1})=0$ for $k<2g$ (see Proposition 2.1 in \cite{BFP}). Then the isomorphism
 $$
H^k(\orgc_{d})=\prod_{g\geq 1} H_c^{k-1+2g(d-1)}(\cM_{g,1})
$$ 
implies that  $H^0(\orgc_{d})=0$ for $-1+2g(d-1)<2g$, i.e. for $d\leq 2$. By  
Theorem {\ref{1: Theorem on H(Der)}} we conclude that
$$
H^0(\Der(\PCY_d)=\K \ \text{for} \ d\leq 2.
$$
This result implies that the dg properad $\PCY_d$ for $d\leq 2$ has no homotopy non-trivial 
automorphisms (except the standard rescaling) which preserve the number of boundaries of the generating ribbon graphs with hairs. 

\sip

There is an  injection \cite{CGP}
$$
\prod_{g\geq 2} H^{\bu}(\cM_{g}) \hook \prod_{g\geq 2} H^{\bu}(\cM_{g,1}) \hook
\prod_{g\geq 2} H^{\bu+2}(\cM_{g,1}) 
$$
where the second map is the cup product with the Euler class. In terms of the compactly supported
cohomology one can say that $\prod_{g\geq 2} H_c^{\bu}(\cM_{g,1})$ contains a subspace isomorphic
to $\prod_{g\geq 2} H_c^{\bu}(\cM_{g})$.  It was proven in \cite{CGP} that the totality of cohomology groups 
$$
\prod_{g\geq 2}H_c^{k+2g}(\cM_g)
$$
 contains a subspace isomorphic to $H^k(\GC_2)$,
where $\GC_2$ is the Kontsevich graph complex (see also \cite{AWZ} for a very short and beautiful proof
of this result). We conclude that for $d=2$ the vector space
$$
H^1(\Der(\PCY_2))= H^1(\orgc_2) \simeq \prod_{g\geq 1} H_c^{2g}(\cM_{g,1})
$$
contains a subspace isomorphic to $H^0(\GC_2)$ which, as was proven in \cite{Wi1}, can be identified with the 
Grothendieck-Teichm\"uller Lie algebra $\grt_1$.

\subsection{Non-Koszulness of the properad $\mathsf{BIB}_d$}

 Consider in $\PCY_d$ the differential closure $\hat{I}$ of an ideal $I$ generated
by haired ribbon graphs having at least one vertex of valency $\geq 4$. The quotient properad 
$$
\mathsf{BIB}_d:= \PCY_d/\hat{I}
$$
has a trivial differential. This properad has been introduced and studied in \cite{Q}  where it was called the properad of balanced infinitesimal bialgebras. It is generated by two sets of planar corollas,
$$
\K[\bS_2][1]=\mbox{span}\left\langle
\Ba{c}\begin{xy}
 <0mm,0.66mm>*{};<0mm,3mm>*{}**@{-},
 <0.39mm,-0.39mm>*{};<2.2mm,-2.2mm>*{}**@{-},
 <-0.35mm,-0.35mm>*{};<-2.2mm,-2.2mm>*{}**@{-},
 <0mm,0mm>*{\circ};<0mm,0mm>*{}**@{},
   <0.39mm,-0.39mm>*{};<2.9mm,-4mm>*{^{_2}}**@{},
   <-0.35mm,-0.35mm>*{};<-2.8mm,-4mm>*{^{_1}}**@{},
\end{xy}\Ea,
\Ba{c}\begin{xy}
 <0mm,0.66mm>*{};<0mm,3mm>*{}**@{-},
 <0.39mm,-0.39mm>*{};<2.2mm,-2.2mm>*{}**@{-},
 <-0.35mm,-0.35mm>*{};<-2.2mm,-2.2mm>*{}**@{-},
 <0mm,0mm>*{\circ};<0mm,0mm>*{}**@{},
   <0.39mm,-0.39mm>*{};<2.9mm,-4mm>*{^{_1}}**@{},
   <-0.35mm,-0.35mm>*{};<-2.8mm,-4mm>*{^{_2}}**@{},
\end{xy}\Ea
\right\rangle,
\ \ \ \
\K[\bS_2][d-2]=\mbox{span}\left\langle
\Ba{c}\begin{xy}
 <0mm,-0.55mm>*{};<0mm,-2.5mm>*{}**@{-},
 <0.5mm,0.5mm>*{};<2.2mm,2.2mm>*{}**@{-},
 <-0.48mm,0.48mm>*{};<-2.2mm,2.2mm>*{}**@{-},
 <0mm,0mm>*{\circ};<0mm,0mm>*{}**@{},
 <0.5mm,0.5mm>*{};<2.7mm,2.8mm>*{^{_2}}**@{},
 <-0.48mm,0.48mm>*{};<-2.7mm,2.8mm>*{^{_1}}**@{},
 \end{xy}\Ea,
\Ba{c}\begin{xy}
 <0mm,-0.55mm>*{};<0mm,-2.5mm>*{}**@{-},
 <0.5mm,0.5mm>*{};<2.2mm,2.2mm>*{}**@{-},
 <-0.48mm,0.48mm>*{};<-2.2mm,2.2mm>*{}**@{-},
 <0mm,0mm>*{\circ};<0mm,0mm>*{}**@{},
 <0.5mm,0.5mm>*{};<2.7mm,2.8mm>*{^{_1}}**@{},
 <-0.48mm,0.48mm>*{};<-2.7mm,2.8mm>*{^{_2}}**@{},
 \end{xy}\Ea
   \right\rangle
$$
modulo a non-commutative version of Drinfeld's relations for a (degree shifted)
Lie bialgebra which are given explicitly in \S 2.1 of \cite{Q}.

\subsubsection{\bf Proposition} {\it The properad $\mathsf{BIB}_d$ is not Koszul}.

\begin{proof} Assume the contrary, i.e.\ assume that the  epimorphism
\Beq\label{5: projection from PCY to BIB}
p: \PCY_d \lon \mathsf{BIB}_d
\Eeq
is a quasi-isomorphism. There is an isomorphism of complexes,
$$
 \Der^{full}(\PCY_d) = \Def^{full}\left(\PCY_d \rar \PCY_d\right)[1],
$$
where the  complex in the r.h.s.\ stands for the standard (i.e. with no ``boundary preserving" restrictions as in \S 5.6)  deformation 
complex of the identity morphism (see \cite{MV}). 
 One can also consider  a deformation complex
$$
\Def^{full}(\PCY_d \lon \mathsf{BIB}_d) 
$$
of the projection $p$ which is generated by at most trivalent haired ribbon graphs. If $p$ is quasi-isomorphism, then the natural epimorphism of complexes,
$$
p^{ind}: \Def^{full}(\PCY_d \rar \PCY_d) \lon \Def^{full}(\PCY_d \lon \mathsf{BIB}_d) 
$$
is also a quasi-isomorphism. The complex $\Def^{full}(\PCY_d \lon \mathsf{BIB}_d)[1]$
is best understood as the complex $\Der^{full}(\PCY_d \lon \mathsf{BIB}_d)$ of {\it derivations of $\PCY_d$ with values in $\mathsf{BIB}_d$}, which is spanned by linear equivariant maps
$$
D:\PCY_d \lon \mathsf{BIB}_d
$$
satisfying the condition
$$
D(\Ga_1\circ \Ga_2)=D(\Ga_1)\circ p(\Ga_2) +(-1)^{|D|} \p(\Ga_1)\circ D(\Ga_2) 
$$
for any properadic composition $\Ga_1\circ \Ga_2$ of any elements $\Ga_1,\Ga_2\in \PCY_d$.
The above epimorphism $p^{ind}$  can be re-written as an epimorphism of the following  complexes,
$$
p^{ind}: \Der^{full}(\PCY_d) \lon \Der^{full}(\PCY_d \lon \mathsf{BIB}_d).
$$
Denoting the image of the subcomplex  $\Der(\PCY_d)\subset  \Der^{full}(\PCY_d)$ under this map
by $\Der(\PCY_d \lon \mathsf{BIB}_d)$, we obtain an epimorphism of complexes
$$
\pi^{ind}: \Der(\PCY_d)\lon \Der(\PCY_d \lon \mathsf{BIB}_d)
$$
which are going to consider in more detail. The main point is that the latter map is also a quasi-isomorphism if the map  (\ref{5: projection from PCY to BIB}) is a quasi-isomorphism
(the argument is the same as in the ``full" case --- use a filtration of both sides by the total number of hairs which picks up on the initial pages the summand $\delta$ of the full differential $d$ in 
(\ref{5: d in Der(PCY)})). This fact prompts us to consider  a quotient complex generated by equivalence classes of at most trivalent ribbon quivers,
$$
\orgc_d^{\mathsf{T}}:= \orgc_d/\widehat{A},
$$
where $\widehat{A}$ stands  the differential closure of the subspace $A\subset \orgc_d$ generated by ribbon quivers having at least one vertex of valency $\geq 4$. As $\pi^{ind}$ is a quasi-isomorphism, we conclude that there is an injection of cohomology groups for any $k$,
\Beq\label{5: injection for H(M)}
\prod_{g\geq 1} H_c^{k-1+2g(d-1)}(\cM_{g,1})=H^k(\orgc_d)\hook H^k(\orgc_d^{\mathsf{T}}),
\Eeq
i.e.\ every cohomology class in the l.h.s.\ can be represented by a ribbon quiver 
$\Ga\in \orgc_d^{\mathsf{T}}$ which has one boundary and at most trivalent vertices.
Assume $\Ga$ has $p_2$ bivalent vertices and $p_3$ trivalent vertices. Then
$$
\# V(\Ga)=p_2+p_3,\ \ \ \#E(\Ga)=\frac{1}{2}(2p_2 + 3p_3)=p_2+ \frac{3}{2}p_3, \ \ \  \#E(\Ga)-\# V(\Ga)=\frac{1}{2}p_3.
$$
As $\Ga$ has one boundary, its genus is given by the formula
$$
 2g= 2+\left(\# E(\Ga) - \# V(\Ga)- \# B(\Ga)\right)=\# E(\Ga) - \# V(\Ga)+ 1= 
 \frac{1}{2}p_3+1
 $$
Hence the cohomological degree of a genus $g$ ribbon quiver 
$\Ga\in \orgc_d^{\mathsf{T}}$ satisfies the inequality,
\Beqrn
|\Ga| &=& d(\#V(\Ga)-1) + (1-d)\#E(\Ga)\\
      &=& d( \# V(\Ga)-\# E(\Ga) -1) + \# E(\Ga) \\
      &=& -2gd + p_2 +  (6g-3)\\
      &\geq & -2gd + 6g-1,
\Eeqrn
where we used the fact that $\Ga$, being oriented, must have at least one bivalent source and at least one bivalent
target, i.e.\ $p_2\geq 2$. The injection (\ref{5: injection for H(M)}) says that every non-zero cohomology class in  $H^p_c(\cM_{g,1})$ can be represented by a genus $g$ ribbon quiver $\Ga\in \orgc_d^{\mathsf{T}}$ of degree
$$
|\Ga| = p+1 -2g(d-1)
$$
which is greater than or equal to $-2gd + 6g-1$, i.e.\ $p$ must satisfy the inequality
$$
p=|\Ga| -1+2g(d-1)\geq -2gd + 6g-1-1+2g(d-1)=4g-2,
$$
implying in turn that $H_c^{p}(\cM_{g,1})$ must vanish for $p<4g-2$.
This conclusion contradicts, e.g., the fact that $H^4_c(\cM_{2,1}) \neq 0$. Hence the assumption that $\mathsf{BIB}_d$ is Koszul is wrong.
\end{proof}

   The ``commutative" version $\OGC_d$ of the dg Lie algebra $\orgc_d$ admits a ``small"
   model generated by equivalence classes of trivalent graphs
   because the properad $\LB_d$ of degree $d$ Lie bialgebras is Koszul  (see \cite{Me2}). 
   The above
   result says that no such ``small" model exists for the dg Lie algebra $\orgc_d$
   of ribbon quivers with one boundary, and hence for the compactly supported cohomology 
   $\prod_{g\geq 1} H_c^{\bu}(\cM_{g,1})$.

\def\cprime{$'$}


\begin{thebibliography}{10}



\bibitem[AWZ]{AWZ} Assar Andersson, Thomas Willwacher and Marko \v Zivkovi\' c,
{\it Oriented hairy graphs and moduli spaces of curves},  arXiv:2005.00439 (2020).

\bibitem[BFP]{BFP} Jonas Bergstr{\"o}m, Carel Faber and Sam Payne, {\it Polynomial point counts and odd cohomology vanishing on
moduli spaces of stable curves}, Ann. Math.  {\bf 199} (2024), 1323-1365.

\bibitem[CGP]{CGP}
Melody Chan, Soren Galatius and Sam Payne, {\it  Tropical curves, graph complexes, and top weight cohomology
of $\cM_g$}, J.\ Amer.\ Math.\ Soc.\ {\bf 34} (2021),  565-594.


\bibitem[C]{Co} Kevin Costello, {\em
A dual version of the ribbon graph decomposition
of moduli space}, Geometry \& Topology {\bf 11} (2007) 1637-1652.

\bibitem[DSV]{DSV}  Vladimir Dotsenko, Sergey Shadrin, and Bruno Vallette, {\em Maurer-Cartan methods in deformation
theory: the twisting procedure},
London Mathematical Society Lecture Note Series (488),
Cambridge University Press, 2023.

\bibitem[IK]{IK} Natalia Iyudu and Maxim Kontsevich, {\em Pre-Calabi-Yau algebras and
noncommutative calculus on higher cyclic Hochschild cohomology},  arXiv:
2011.11888 (2020).

\bibitem[IKV]{IKV} Natalia Iyudu, Maxim Kontsevich, and Yannis Vlassopoulos. {\em Pre-CalabiYau
algebras as noncommutative Poisson structures},  Journal of Algebra {\bf 567} (2021), pp. 63-90.

\bibitem[K1]{Ko1} Maxim Kontsevich, {\em Intersection Theory on the Moduli Space of Curves
and the Matrix Airy Function}, Commun.\ Math.\ Phys. {\bf 147} (1992)  1-23.


 \bibitem[K2]{Ko2} Maxim Kontsevich, {\it Formality Conjecture}, D. Sternheimer et al. (eds.),
Deformation Theory and Symplectic
Geometry, Kluwer 1997, 139-156.

 \bibitem[K3]{Ko3} Maxim Kontsevich, {\it Derived Grothendieck-Teichm\"uller group and 
 graph complexes [after T. Willwacher]}. Seminaire Bourbaki, vol. 2016/2017, 1120-1135. 
 Asterisque No. 407 (2019), Exp. No. 1126, 183-211


\bibitem[KTV]{KTV}
Maxim Kontsevich, Alex Takeda, and Yiannis Vlassopoulos {\it  Pre-Calabi-Yau algebras
and topological quantum field theories},  preprint arXiv: 2112.14667 (2021).


\bibitem[vL]{L}
Pepijn van der Laan, {\em Operads up to Homotopy and
Deformations of Operad Maps}\, arXiv:math.QA/0208041 (2005)


\bibitem[LV]{LV} Johan Leray and Bruno Vallette, {\em Pre-Calabi-Yau algebras and homotopy
double Poisson gebras}, arXiv: 2203.05062 (2022).

\bibitem[M1]{Me1} Sergei Merkulov, {\it Gravity prop and moduli spaces $\cM_{g,n}$},
arXiv:2108.10644 (2021)


\bibitem[M2]{Me2} Sergei Merkulov, {\it On interrelations between graph complexes},
Int.\ Math.\ Res.\ Notices,  {\bf  2025}, Issue 2  (2025), rnae287,
https://doi.org/10.1093/imrn/rnae287.

\bibitem[M3]{Me3} Sergei Merkulov,  {\it Pre-Calabi-Yau algebras and oriented gravity properad},
arXiv:2501.11158 (2025)


\bibitem[MV]{MV}  Sergei Merkulov and  Bruno Vallette,
{\it Deformation theory of representations of prop(erad)s I \& II},
{ Journal f\"ur die reine und angewandte Mathematik (Qrelle)  634}, 51-106,
 \& {\bf 636}, 123-174 (2009)


 \bibitem[MW1]{MW1} Sergei Merkulov and Thomas Willwacher, {\it Props of ribbon graphs,
 involutive Lie bialgebras and moduli spaces of curves}, preprint  arXiv:1511.07808
  (2015) 51pp.


\bibitem[MW2]{MW2} Sergei Merkulov and Thomas Willwacher, {\it Deformation theory of
 Lie bialgebra properads},   In: Geometry and Physics: A Festschrift in honour of
 Nigel Hitchin, Oxford University Press 2018, pp. 219-248




\bibitem[MWW]{MWW} Sergei Merkulov, Thomas Willwacher and Vincent Wolff,
{\it The oriented graph complex revisited}, arXiv:2411.19657 (2024)


\bibitem[P]{Pe}  Robert Penner, {\it The decorated Teichm\"uller space of punctured surfaces},
Comm.\ Math.\
Phys. {\bf 113} (1987), 299-339.



\bibitem[Q]{Q} Alexandre Quesney, {\it Balanced infinitesimal bialgebras, 
double Poisson gebras and pre-Calabi-Yau algebras}, arXiv:2312.14893
(2023)




 \bibitem[W1]{Wi1} Thomas Willwacher, {\it M.\ Kontsevich's graph complex and the
  Grothendieck-Teichmueller Lie algebra}, Invent. Math. {\bf 200} (2015), 671-760.

\bibitem[W2]{Wi2} Thomas Willwacher, {\it Oriented graph complexes},
Comm. Math. Phys. {\bf 334} (2015), no. 3, 1649-1666.

\bibitem[W3]{Wi3} Thomas Willwacher, {\it Little disks operads and Feynman diagrams}. 
Proceedings of the International Congress of Mathematicians - Rio de Janeiro 2018. Vol. 
II. Invited lectures, 1241-1261, World Sci. Publ., Hackensack, NJ, 2018.


\bibitem[Z]{Z} M.\ \v Zivkovi\' c, {\it Multi-oriented graph complexes and quasi-isomorphisms between them I: oriented graphs},
 High.\ Struct. {\bf 4}, no.1 (2020) 266-283.

 \end{thebibliography}
\end{document}